\documentclass[12pt]{article}
\usepackage{subeqn,amsfonts,amssymb}
\newcommand{\CC}{{\mathbb C}}
\newcommand{\HH}{{\mathbb H}}
\newcommand{\II}{{\mathbb I}}
\newcommand{\NN}{{\mathbb N}}
\newcommand{\ZZ}{{\mathbb Z}}
\newcommand{\eps}{{\varepsilon}}
\newcommand{\tht}{{\vartheta}}
\newcommand{\sfrac}[2]{{\textstyle{\frac{#1}{#2}}}}
\newcommand{\car}[2]{\left[\begin{array}{c}{#1}\\{#2}\end{array}\right]}
\newcommand{\finproof}{{\hfill \rule{5pt}{5pt}}}

\newtheorem{thm}{Theorem}
\newtheorem{prop}{Proposition}

\begin{document}
\newpage
\pagestyle{empty}
\setcounter{page}{0}
\vfill
\begin{center}

{\Large {\bf Deformed ${\cal W}_{N}$ algebras from elliptic $sl(N)$
algebras}}

\vspace{10mm}

{\large J. Avan}

\vspace{4mm}

{\em LPTHE, CNRS-URA 280, Universit{\'e}s Paris VI/VII, France}

\vspace{7mm}

{\large L. Frappat, M. Rossi, P. Sorba}

\vspace{4mm}

{\em Laboratoire de Physique Th{\'e}orique ENSLAPP, CNRS-URA 1436,}
\\
{\em {\'E}cole Normale Sup{\'e}rieure de Lyon and Universit{\'e} de Savoie,
France}
\end{center}

\vfill

\begin{abstract}
We extend to the $sl(N)$ case the results that we previously obtained
on the construction of ${\cal W}_{q,p}$ algebras from the
elliptic algebra ${\cal A}_{q,p}(\widehat{sl}(2)_{c})$.  The elliptic
algebra ${\cal A}_{q,p}(\widehat{sl}(N)_{c})$ at the critical level
$c=-N$ has an extended center containing trace-like operators
$t(z)$.  Families of Poisson structures indexed by $N(N-1)/2$ integers,
defining $q$-deformations of the ${\cal W}_{N}$ algebra, are
constructed.  The operators $t(z)$ also close an exchange algebra when
$(-p^\sfrac{1}{2})^{NM} = q^{-c-N}$ for $M\in\ZZ$.  It becomes Abelian
when in addition $p=q^{Nh}$ where $h$ is a non-zero integer.  The
Poisson structures obtained in these classical limits contain
different $q$-deformed ${\cal W}_{N}$ algebras depending on the parity
of $h$, characterizing the exchange structures at $p \ne q^{Nh}$ as new 
${\cal W}_{q,p}(sl(N))$ algebras.
\end{abstract}

\def\abstractname{R{\'e}sum{\'e}}
\begin{abstract}
Nous {\'e}tendons au cas $sl(N)$ les r{\'e}sultats que nous avons obtenus
pr{\'e}c{\'e}demment concernant la construction des alg{\`e}bres ${\cal
W}_{q,p}$ {\`a} partir de l'alg{\`e}bre elliptique ${\cal
A}_{q,p}(\widehat{sl}(2)_{c})$.  L'alg{\`e}bre elliptique ${\cal
A}_{q,p}(\widehat{sl}(N)_{c})$ au niveau critique $c=-N$ poss{\`e}de un
centre \'etendu contenant des op{\'e}rateurs de trace $t(z)$.
On construit sur ce centre des familles de structures de Poisson
indic{\'e}es par $N(N-1)/2$ entiers, d\'efinissant des $q$-d{\'e}formations
de l'alg{\`e}bre ${\cal W}_{N}$.  Les op{\'e}rateurs $t(z)$ engendrent une
alg{\`e}bre d'{\'e}change lorsque $(-p^\sfrac{1}{2})^{NM} = q^{-c-N}$ o{\`u}
$M\in\ZZ$.  Cette alg{\`e}bre devient ab{\'e}lienne si de plus $p=q^{Nh}$
avec $h$ entier non nul.  Les structures de Poisson obtenues dans ces
limites classiques contiennent diff{\'e}rentes alg{\`e}bres ${\cal W}_{N}$
$q$-d{\'e}form{\'e}es d{\'e}pendant de la parit{\'e} de l'entier $h$,
caract{\'e}risant les structures d'\'echange {\`a} $p \ne q^{Nh}$ comme de 
nouvelles alg{\`e}bres ${\cal W}_{q,p}(sl(N))$.
\end{abstract}

\vfill
\vfill

\rightline{ENSLAPP-AL-670/97}
\rightline{PAR-LPTHE 97-49}
\rightline{math.QA/9801105}
\rightline{December 1997}

\newpage
\pagestyle{plain}

\section{Introduction}

The elliptic algebra ${\cal A}_{q,p}(\widehat{sl}(2)_{c})$ was 
introduced in \cite{FIJKMY} and further studied in \cite{JKM,KLP,Konno} 
where several trigonometric limits were derived and shown to be relevant 
as symmetries for the XXZ model, whereas ${\cal 
A}_{q,p}(\widehat{sl}(2)_{c})$ was proposed as basic symmetry algebra 
for the XYZ model \cite{JMN,JKKMW}.  More recently another, possibly 
related elliptic algebra ${\cal U}_{q,p}(\widehat{sl}(2)_{c})$ was 
introduced \cite{Konno2} in connection with the $k$-fusion RSOS model 
\cite{ABF,DJMO}.  A dynamical version of the algebra ${\cal 
A}_{q,p}(\widehat{sl}(2)_{c})$, making use of the dynamical 
``eight-vertex'' elliptic $R$-matrix \cite{ABB96,Fe94}, was also 
recently proposed in \cite{HouYang}.  The domain was unified recently 
\cite{Fron97,JKOS97} where both ${\cal A}_{q,p}(\widehat{sl}(2)_{c})$ 
algebra (and its extension to $sl(N)$) and Felder's elliptic dynamical 
algebra \cite{Fe94,GN84} were shown to be obtained from application of 
universal twisting operators (see also \cite{ABRR97}) to the quantum 
algebra ${\cal U}_{q}(\widehat{sl}(2)_{c})$.  As a consequence the 
construction \cite{FIJKMY} was validated and extended to the $sl(N)$ 
case.

The $q$-Virasoro algebra was introduced in \cite{SKAO} as an extension 
to the Ruijsenaar--Schneider model of the Virasoro algebra arising in 
the collective theory formulation of the Calogero--Moser model 
\cite{AJL,MP,AJ}.  It arose simultaneously as quantization of a 
classical quadratic Poisson structure on the center at $c=-2$ of the 
quantum affine algebra ${\cal U}_{q}(sl(2)_{c})$ \cite{FR}.  It was shown 
\cite{L,LP} to be the symmetry algebra of the restricted ABF model, 
itself connected as previously said to the elliptic ${\cal 
U}_{q,p}(\widehat{sl}(2)_{c})$ algebra.  On the other hand, the limit of 
$q$-Virasoro current coincides \cite{L,JKM} with the field obtained by 
concatenation of vertex operators of the degeneracy limit of ${\cal 
A}_{q,p}(\widehat{sl}(2)_{c})$, hinting indeed to a deep connection 
between the two structures.  This identification was recently extended to 
the full $q$-deformed Virasoro case \cite{JS97}.

A direct connection was indeed established in two recent papers 
\cite{AFRS97a,AFRS97b}: we derived exchange algebras and limit Poisson 
structures from the original ${\cal A}_{q,p}(\widehat{sl}(2)_{c})$ 
algebra, by establishing the existence of a center at $p^m=q^{c+2}$ for 
any integer $m$.  The Poisson structures all boil down to the 
$q$-Virasoro (classical) algebra defined in \cite{SKAO,FR}.  The 
exchange algebras are however quite distinct from the quantization of 
$q$-Virasoro constructed in \cite{SKAO,FF}.  Their eventual central (and 
linear ?)  extensions also exhibit a much richer structure than was 
originally derived in \cite{FF} and are currently under investigation 
\cite{AFRS97c}.

It is known that $q$-Virasoro structures admit extensions to the $sl(N)$ 
case as classical and quantum $q$-${\cal W}_{N}$ algebras 
\cite{FR,FF,AKOS}.  These have been the object of many investigations 
recently \cite{FJMOP,FR2}.  In particular, they were also shown to be in 
general symmetry algebras of RSOS models \cite{AJMP,LP}.  It is a 
natural question in view of our previous results to investigate whether 
classical and quantized $q$-${\cal W}_{N}$ algebras arise as structures 
embedded in the generalized elliptic algebra ${\cal 
A}_{q,p}(\widehat{sl}(N)_{c})$.

\medskip

The plan of the paper is as follows.  We first recall in Sect.  
\ref{sect2} useful properties of elliptic functions and elliptic $sl(N)$ 
$R$-matrix, and introduce the definition of the quantum elliptic algebra 
${\cal A}_{q,p}(\widehat{sl}(N)_{c})$ \cite{JKOS97}.  We then prove in Sect.  
\ref{sect3} the existence of an extended center at $c=-N$ and compute 
the Poisson structures on this center, using the notion of 
sector-depending mode Poisson bracket already introduced in 
\cite{AFRS97a}.  These structures are identified as $q$-deformed ${\cal 
W}_{N}$ algebras.  They are obtained from analytic continuations of a 
classical algebra identical to the initial version of $q$-${\cal W}_{N}$ 
given in \cite{FR}.  In Sect.  \ref{sect4} we show the existence of 
closed (quadratic) exchange algebras whenever $(-p^\sfrac{1}{2})^{NM} = 
q^{-c-N}$ for any integer $M\in\ZZ$.  These algebras differ from the 
quantum ${\cal W}_{q,p}(\widehat{sl}(N))$ structures introduced in 
\cite{FF}.  They admit a classical limit (commuting algebras) at $p = 
q^{Nh}$ with $h \in \ZZ \backslash \{ 0 \}$.  In Sect.  \ref{sect5} we 
compute the related Poisson structures.  They include for $h$ even the 
Poisson structures in \cite{FR}.  The exchange algebras therefore 
realize new quantizations of these Poisson structures.  When $h$ is odd, 
by contrast, this classical limit takes a form different from the 
initial $q$-${\cal W}_{N}$ structures.  This emphasizes the key role of 
the initial 3-parameter structure ${\cal A}_{q,p}(\widehat{sl}(N)_{c})$ 
in allowing for an intermediate quantum $q$-deformed ${\cal W}_{N}$ 
algebra.  Finally in Sect.  \ref{sect6}, we compute the mode expansion 
of the quantum exchange algebra structures.  As in 
the classical case, a ``sector-type'' labeling is needed due to the 
singularity behaviour of the structure function viewed as an analytic 
continuation and therefore exhibiting different formal series 
expansions corresponding to different convergent series expansions in 
distinct domains of the complex plane.  We give an explicit example of 
this treatment applied to the spin one field for the sake of simplicity, 
and we describe the essential features of the extension of our computation 
to higher spin fields.

\section{Notations and basic definitions}
\label{sect2}
\setcounter{equation}{0}

\subsection{Definition of the $N$-elliptic $R$-matrix}

The $N$-elliptic $R$-matrix in $\mbox{End}(\CC^N) \otimes
\mbox{End}(\CC^N)$, associated to the $\ZZ_{N}$-vertex model, is defined as 
follows \cite{Bela,ChCh}:
\begin{equation}
\widetilde R(z,q,p) = z^{2/N-2} \frac{1}{\kappa(z^2)} 
\frac{\tht\car{\sfrac{1}{2}}{\sfrac{1}{2}}(\zeta,\tau)} 
{\tht\car{\sfrac{1}{2}}{\sfrac{1}{2}}(\xi+\zeta,\tau)} \,\, 
\sum_{(\alpha_1,\alpha_2)\in\ZZ_N\times\ZZ_N} 
W_{(\alpha_1,\alpha_2)}(\xi,\zeta,\tau) \,\, I_{(\alpha_1,\alpha_2)} 
\otimes I_{(\alpha_1,\alpha_2)}^{-1} \,,
\label{eq27}
\end{equation}
where the variables $z,q,p$ are related to the variables 
$\xi,\zeta,\tau$ by
\begin{equation}
z=e^{i\pi\xi} \,,\qquad q=e^{i\pi\zeta} \,,\qquad p=e^{2i\pi\tau} \,.
\label{eq28}
\end{equation}
The Jacobi theta functions with rational characteristics
$\tht\car{\gamma_1}{\gamma_2}(\xi,\tau)$ are defined in Appendix A.
\\
The normalization factor is chosen as follows:
\begin{equation}
\frac{1}{\kappa(z^2)} = \frac{(q^{2N}z^{-2};p,q^{2N})_\infty
\, (q^2z^2;p,q^{2N})_\infty \, (pz^{-2};p,q^{2N})_\infty \,
(pq^{2N-2}z^2;p,q^{2N})_\infty} {(q^{2N}z^2;p,q^{2N})_\infty
\, (q^2z^{-2};p,q^{2N})_\infty \, (pz^2;p,q^{2N})_\infty \,
(pq^{2N-2}z^{-2};p,q^{2N})_\infty} \,.
\label{eq29}
\end{equation}
The functions $W_{(\alpha_1,\alpha_2)}$ are given by
\begin{equation}
W_{(\alpha_1,\alpha_2)}(\xi,\zeta,\tau) =
\frac{\tht\car{\sfrac{1}{2}+\alpha_1/N}
{\sfrac{1}{2}+\alpha_2/N}(\xi+\zeta/N,\tau)}
{N\tht\car{\sfrac{1}{2}+\alpha_1/N}
{\sfrac{1}{2}+\alpha_2/N}(\zeta/N,\tau)} \,.
\label{eq211}
\end{equation}
The matrices $I_{(\alpha_1,\alpha_2)}$ are defined as follows:
\begin{equation}
I_{(\alpha_1,\alpha_2)} = g^{\alpha_2} \, h^{\alpha_1} \,,
\label{212}
\end{equation}
where the $N \times N$ matrices $g$ and $h$ are given by $g_{ij} =
\omega^i\delta_{ij}$ and $h_{ij} = \delta_{i+1,j}$, the addition of 
indices being understood modulo $N$.  \\
Let us set
\begin{equation}
S(\xi,\zeta,\tau) = \sum_{(\alpha_1,\alpha_2)\in\ZZ_N\times\ZZ_N}
W_{(\alpha_1,\alpha_2)}(\xi,\zeta,\tau) \,\, I_{(\alpha_1,\alpha_2)}
\otimes I_{(\alpha_1,\alpha_2)}^{-1} \,.
\label{eq210}
\end{equation}
The matrix $S$ is $\ZZ_N$-symmetric, that is
\begin{equation}
S_{c+s\,,\,d+s}^{a+s\,,\,b+s} = S_{c\,,\,d}^{a\,,\,b}
\label{eq213}
\end{equation}
for any indices $a,b,c,d,s \in \ZZ_N$ (the addition of indices being
understood modulo $N$) and the non-vanishing elements of the matrix
$S$ are $S_{c\,,\,a+b}^{a\,,\,c+b}$.  One finds explicitly:
\begin{equation}
S_{c\,,\,a+b}^{a\,,\,c+b} = \sum_{s\in\ZZ_N} W_{(a-c,s)} \,
\omega^{-bs} \,.
\label{eq214}
\end{equation}
The matrix $S$ being $\ZZ_N$-symmetric, it is sufficient to examine the
terms $S^{ab} \equiv S^{0\,,\,a+b}_{a\,,\,b}$.
One finds:
\begin{equation}
S^{ab}(\xi,\zeta,\tau) =
\frac{\tht\car{(b-a)/N+\sfrac{1}{2}}{\sfrac{1}{2}}(\xi+\zeta,N\tau)}
{\tht\car{-a/N+\sfrac{1}{2}}{\sfrac{1}{2}}(\zeta,N\tau)} \,\,
\frac{\displaystyle\prod_{k=0,k \ne b}^{N-1}
\tht\car{k/N+\sfrac{1}{2}}{\sfrac{1}{2}}(\xi,N\tau)}
{\displaystyle\prod_{k=1}^{N-1}\tht\car{k/N+\sfrac{1}{2}}
{\sfrac{1}{2}}(0,N\tau)} \,.
\label{eq215}
\end{equation}
It satisfies the following shift properties:
\begin{subequations}
\begin{eqnarray}
&& S^{ab}(\xi,\zeta+\lambda\tau,\tau) = \exp(-2i\pi\lambda \xi/N)
\, S^{a-\lambda,b}(\xi,\zeta,\tau) \,, \\
&& S^{ab}(\xi+\lambda\tau,\zeta,\tau) =
\exp(-i\pi\lambda^2\tau-2i\pi\lambda (\xi+\zeta/N+\sfrac{1}{2})) \,
S^{a,b+\lambda}(\xi,\zeta,\tau) \,.
\end{eqnarray}
\label{eq216}
\end{subequations}
Using the following ``gluing'' formula:
\begin{equation}
\prod_{k=0}^{N-1}\tht\car{k/N+\sfrac{1}{2}}{\sfrac{1}{2}}(\xi,N\tau)
= p^{\frac{1}{24}(N^2-1)} \frac{\prod_{k=0}^{N-1}
(p^{N-k};p^N)_\infty}{(p^N;p^N)_\infty^N} \,\,
\tht\car{\sfrac{1}{2}}{\sfrac{1}{2}}(\xi,\tau) \,,
\label{eq217}
\end{equation}
and inserting eq.  (\ref{eq215}) into eqs.  (\ref{eq210}) and
(\ref{eq27}), one then finds the following expression in terms of the
Jacobi $\Theta$ functions for the $\widetilde R^{ab}$ elements of the
matrix (\ref{eq27}):
\begin{equation}
\widetilde R^{ab}(z,q,p) = \frac{1}{\kappa(z^2)} \,\, p^{-\frac{ab}{N}} 
q^{\frac{2b}{N}} z^{-\frac{2}{N}(N+a-1)} \,\, 
\frac{\Theta_{p^N}(p^{N+b-a}q^2z^2)\,\Theta_{p^N}(p^N)} 
{\Theta_{p^N}(p^{N+b}z^2)\,\Theta_{p^N}(p^{N-a}q^2)} 
\frac{\Theta_p(pq^2)\,\Theta_p(pz^2)}{\Theta_p(pq^2z^2)\,\Theta_p(p)} \,.
\label{eq219}
\end{equation}

\subsection{Gauge-transformed $R$-matrix}

In order to make the comparison with our previous results easier
\cite{AFRS97a,AFRS97b}, one needs to introduce the following
``gauge-transformed'' matrix:
\begin{equation}
R(z,q,p) = (g^{\frac 12} \otimes g^{\frac 12}) \widetilde
R(z,q,p) (g^{-\frac 12} \otimes g^{-\frac 12}) \,.
\label{eq220}
\end{equation}
\begin{thm}\label{thmone}
The matrix $R(z,q,p)$ satisfies the following properties:
\\
-- Yang--Baxter equation:
\begin{equation}
R_{12}(z) \, R_{13}(w) \, R_{23}(w/z) =
R_{23}(w/z) \, R_{13}(w) \, R_{12}(z) \,,
\label{eq221}
\end{equation}
-- Unitarity:
\begin{equation}
R_{12}(z) \, R_{21}(z^{-1}) = 1 \,,
\label{eq222}
\end{equation}
-- Crossing-symmetry:
\begin{equation}
R_{12}(z)^{t_2} \, R_{21}(z^{-1}q^{-N})^{t_2} = 1 \,,
\label{eq223}
\end{equation}
-- Antisymmetry:
\begin{equation}
R_{12}(-z) = \omega \, (g^{-1} \otimes \II) \, R_{12}(z) \,
(g \otimes \II) \,,
\label{eq224}
\end{equation}
-- Quasi-periodicity:
\begin{equation}
\widehat R_{12}(-z p^{\frac 12}) = (g^{\frac 12} h g^{\frac 12}
\otimes \II)^{-1} \, \widehat R_{21}(z^{-1})^{-1} \,
(g^{\frac 12} h g^{\frac 12} \otimes \II) \,,
\label{eq225}
\end{equation}
where
\begin{equation}
\widehat R_{12}(z) \equiv \widehat R_{12}(z,q,p) =
\tau_N(q^{\frac 12}z^{-1}) \, R_{12}(z,q,p) \,,
\label{eq226}
\end{equation}
the function $\tau_N(z)$ being defined by
\begin{equation}
\tau_N(z) = z^{\frac{2}{N}-2} \,
\frac{\Theta_{q^{2N}}(qz^2)}{\Theta_{q^{2N}}(qz^{-2})} \,.
\label{227}
\end{equation}
The function $\tau_N(z)$ is periodic with period $q^N$:
$\tau_N(q^Nz) = \tau_N(z)$ and satisfies $\tau_N(z^{-1}) =
\tau_N(z)^{-1}$.
\end{thm}

\medskip

\noindent
{\bf Proof:} The proof of the Yang--Baxter equation has been given
in \cite{Tra85}.  The proof of the unitarity and the
crossing-symmetry is done by a direct calculation.  One has to use
the following two identities (the first one for the unitarity and
the second one for the crossing-symmetry, see \cite{RT86} for
a proof of eqs.  (\ref{eq228}) and (\ref{eq229})):
\begin{eqnarray}
&& \hspace{-10mm} \sum_{k\in\ZZ_N} S^{-i-k,i-k}(\xi,\zeta,\tau)
S^{-j-k,k-j}(-\xi,\zeta,\tau) = \frac{\tht\car{\sfrac{1}{2}}
{\sfrac{1}{2}}(\xi+\zeta,\tau)\,\tht\car{\sfrac{1}{2}}
{\sfrac{1}{2}}(-\xi+\zeta,\tau)}
{\tht\car{\sfrac{1}{2}}{\sfrac{1}{2}}(\zeta,\tau)^2} \,
\delta_{ij} \,, \label{eq228} \\
&& \hspace{-10mm} \sum_{k\in\ZZ_N} S^{i-k,0}(z,\zeta,\tau)
S^{j-k,0}(-\xi-N\zeta,\zeta,\tau) = \frac{\tht\car{\sfrac{1}{2}}
{\sfrac{1}{2}}(z,\tau)\,\tht\car{\sfrac{1}{2}}
{\sfrac{1}{2}}(-\xi-N\zeta,\tau)}
{\tht\car{\sfrac{1}{2}}{\sfrac{1}{2}}(\zeta,\tau)^2} \,
\delta_{ij} \,. \label{eq229}
\end{eqnarray}
Finally, the antisymmetry and the quasi-periodicity are explicitly
checked from the expressions of the matrix elements of $R$.
\finproof

\medskip

\noindent
{\bf Remark}: The crossing-symmetry and the unitarity properties of
$R_{12}$ allow exchange of inversion and transposition for the
matrix $R_{12}$ as (the same property also holds for the matrix
$\widehat R_{12}$):
\begin{equation}
\Big(R_{12}(x)^{t_2}\Big)^{-1} = \Big(R_{12}(q^Nx)^{-1}\Big)^{t_2} \,.
\label{eq230}
\end{equation}

\subsection{The quantum elliptic algebra ${\cal A}_{q,p}(\widehat{sl}(N)_{c})$}

We now define the elliptic quantum algebra ${\cal 
A}_{q,p}(\widehat{sl}(N)_{c})$ \cite{FIJKMY,JKOS97} as an algebra of 
operators $L_{ij}(z)$ $\equiv \sum_{n\in\ZZ} L_{ij}(n) \, z^n$ 
where $i,j \in \ZZ_N$, encapsulated 
into a $N \times N$ matrix
\begin{equation}
L(z) = \left(\begin{array}{ccc} L_{11}(z) & \cdots &
L_{1N}(z) \cr \vdots && \vdots \cr L_{N1}(z) & \cdots &
L_{NN}(z) \cr \end{array}\right) \,.
\label{eq231}
\end{equation}
One defines ${\cal A}_{q,p}(\widehat{gl}(N)_c)$ by imposing the
following constraints on the $L_{ij}(z)$ (with the matrix
$\widehat R_{12}$ given by eq. (\ref{eq226})):
\begin{equation}
\widehat R_{12}(z/w) \, L_1(z) \, L_2(w) =
L_2(w) \, L_1(z) \, \widehat R_{12}^{*}(z/w) \,,
\label{eq232}
\end{equation}
where $L_1(z) \equiv L(z) \otimes \II$, $L_2(z) \equiv \II
\otimes L(z)$ and $\widehat R^{*}_{12}$ is defined by $\widehat
R^{*}_{12}(z,q,p) \equiv \widehat R_{12}(z,q,p^*=pq^{-2c})$.  
This definition is the most immediate generalization to $N$ of the 
definition adopted in \cite{FIJKMY} for $N=2$.  \\
The matrix $\widehat R^{*}_{12}$ obeys also the unitarity,
crossing-symmetry, antisymmetry and quasi-periodicity conditions of
Theorem \ref{thmone} (note that the quasi-periodicity condition
(\ref{eq225}) for $\widehat R^{*}_{12}$ has to be understood with
the modified elliptic nome $p^*$).  \\
The $q$-determinant $q$-$\det L(z)$ given by
\begin{equation}
q\mbox{-}\det L(z) \equiv \sum_{\sigma\in\mathfrak S_N}
\eps(\sigma) \prod_{i=1}^N L_{i,\sigma(i)}(z q^{i-N-1})
\label{eq233}
\end{equation}
($\eps(\sigma)$ being the signature of the permutation $\sigma$) is in
the center of ${\cal A}_{q,p}(\widehat{gl}(N)_c)$. It can be
``factored out'', and set to the value $q^{\frac c2}$ so as to get
\begin{equation}
{\cal A}_{q,p}(\widehat{sl}(N)_c) = {\cal A}_{q,p}(\widehat{gl}(N)_c)/
\langle q\mbox{-}\det L - q^{\frac c2} \rangle \,.
\label{eq234}
\end{equation}
It is useful to introduce the following two matrices:
\begin{subequations}
\begin{eqnarray}
&& L^+(z) \equiv L(q^{\frac c2}z) \,, \\
&& L^-(z) \equiv (g^{\frac 12} h g^{\frac 12}) \,
L(-p^{\frac 12}z) \, (g^{\frac 12} h g^{\frac 12})^{-1} \,.
\end{eqnarray}
\label{eq235}
\end{subequations}
They obey coupled exchange relations following from (\ref{eq232}) and
periodicity/unitarity properties of the matrices
$\widehat R_{12}$ and $\widehat R^{*}_{12}$:
\begin{subequations}
\begin{eqnarray}
&& \widehat R_{12}(z/w) \, L^\pm_1(z) \,
L^\pm_2(w) = L^\pm_2(w) \, L^\pm_1(z) \,
\widehat R^{*}_{12}(z/w) \,, \\
&& \widehat R_{12}(q^{\frac c2}z/w) \,
L^+_1(z) \, L^-_2(w) = L^-_2(w) \, L^+_1(z)
\, \widehat R^{*}_{12}(q^{-\frac c2}z/w) \,.
\end{eqnarray}
\label{eq236}
\end{subequations}
The parameters $c,p,q$ in our definition are related to the corresponding 
parameters $c',p',q'$ of \cite{JKOS97} by $c=2c'$, $p={(p')}^{2/N}$, 
$q={(q')}^{1/N}$.

\section{The center of ${\cal A}_{q,p}(\widehat{sl}(N)_{c})$ at the
critical level $c=-N$}
\label{sect3}
\setcounter{equation}{0}

\subsection{Center of ${\cal A}_{q,p}(\widehat{sl}(N)_{c})$}

\begin{thm}\label{thmtwo}
At the critical level $c=-N$, the operators generated by
\begin{equation}
t(z) = {\rm Tr}(L(z)) = {\rm Tr}\Big(L^+(q^{\frac c2}z)
L^-(z)^{-1}\Big)
\label{eq31}
\end{equation}
lie in the center of the algebra ${\cal A}_{q,p}(\widehat{sl}(N)_{c})$.
\end{thm}

\medskip

\noindent
{\bf Proof:}
Defining $\widetilde L^\pm(z) \equiv (L^\pm(z)^{-1})^t$, one can
derive from eqs. (\ref{eq236}) further exchange relations between the
operators $L^+$ and $\widetilde L^-$:
\begin{subequations}
\begin{eqnarray}
&& \left(\widehat R_{12}(z/w)^{t_2}\right)^{-1} \, L_1^\pm(z) \,
\widetilde L_2^\pm(w) = \widetilde L_2^\pm(w) \, L_1^\pm(z) \,
\left({\widehat R_{12}^*(z/w)}^{t_2}\right)^{-1} \,,
\label{eq32a} \\
&& \left(\widehat R_{12}(q^{\frac c2}z/w)^{t_2}\right)^{-1} \,
L_1^+(z) \, \widetilde L_2^-(w) = \widetilde L_2^-(w) \, L_1^+(z) \,
\left(\widehat R_{12}^*(q^{-\frac c2}z/w)^{t_2}\right)^{-1} \,.
\label{eq32b}
\end{eqnarray}
\label{eq32}
\end{subequations}
Let us now compute $\Big[ t(z),L^+(w) \Big]$. One rewrites:
\begin{equation}
t(z) \, L^+(w) = {\rm Tr}_1 \Big(L_1^+(q^{\frac c2}z) \widetilde
L_1^-(z)^{t_1}\Big) L_2^+(w) = {\rm Tr}_1
\Big(L_1^+(q^{\frac c2}z)^{t_1} \widetilde L_1^-(z) L_2^+(w)\Big)
\label{eq33}
\end{equation}
since one is allowed to exchange transposition under a trace procedure.
Commuting $L_2^+(w)$ through $\widetilde L_1^-(z)$ using eq.
(\ref{eq32b}), one gets:
\begin{equation}
t(z) \, L_2^+(w) = {\rm Tr}_1\Big(L_1^+(q^{\frac c2}z)^{t_1}
(\widehat R_{21}(q^{\frac c2}w/z)^{t_1})^{-1} L_2^+(w) \widetilde
L_1^-(z) \widehat R_{21}^*(q^{-\frac c2}w/z)^{t_1}\Big) \,.
\label{eq34}
\end{equation}
Using the unitarity property of $\widehat R_{12}$, one has:
\begin{equation}
L_1^+(q^{\frac c2}z)^{t_1} \,
(\widehat R_{21}(q^{-\frac c2}w/z)^{-1})^{t_1} \, L_2^+(w) =
L_2^+(w) \, (\widehat R_{21}^*(q^{-\frac c2}w/z)^{-1})^{t_1}
\, L_1^+(q^{\frac c2}z)^{t_1} \,.
\label{eq35}
\end{equation}
Then applying eq. (\ref{eq230}) to (\ref{eq35}) gives:
\begin{equation}
L_1^+(q^{\frac c2}z)^{t_1} \,
(\widehat R_{21}(wq^{-\frac c2 - N}/z)^{t_1})^{-1} \, L_2^+(w) =
L_2^+(w) \, (\widehat R_{21}^*(q^{-\frac c2}w/z)^{-1})^{t_1} \,
L_1^+(q^{\frac c2}z)^{t_1} \,,
\label{eq36}
\end{equation}
which leads at the critical value $c=-N$ to
\begin{equation}
L_1^+(q^{\frac c2}z)^{t_1} \,
(\widehat R_{21}(q^{\frac c2}w/z)^{t_1})^{-1} \, L_2^+(w) =
L_2^+(w) \, (\widehat R_{21}^*(q^{-\frac c2}w/z)^{-1})^{t_1} \,
L_1^+(q^{\frac c2}z)^{t_1} \,.
\label{eq37}
\end{equation}
Now, inserting (\ref{eq37}) into (\ref{eq34}), the first three
terms in (\ref{eq34}) can be rearranged and one obtains:
\begin{equation}
t(z) \, L_2^+(w) = L_2^+(w) \, {\rm Tr}_1\Big(
(\widehat R_{21}^*(q^{-\frac c2}w/z)^{-1})^{t_1}
L_1^+(q^{\frac c2}z)^{t_1} \widetilde L_1^-(z)
\widehat R_{21}^*(q^{-\frac c2}w/z)^{t_1} \Big)
\label{eq38}
\end{equation}
Using the fact that under a trace over the space 1 one has
${\rm Tr}_1 \Big( R_{21} Q_1 {R'}_{21} \Big) = {\rm Tr}_1
\Big( Q_1 {R'_{21}}^{t_2} {R_{21}}^{t_2} \Big)^{t_2}$, one gets
\begin{equation}
t(z) \, L_2^+(w) = L_2^+(w) \, {\rm Tr}_1\Big(
L_1^+(q^{\frac c2}z)^{t_1} \widetilde L_1^-(z)
\widehat R_{21}^*(q^{-\frac c2}w/z)^{t_1t_2}
(\widehat R_{21}^*(q^{-\frac c2}w/z)^{-1})^{t_1t_2} \Big)^{t_2}
\label{eq39}
\end{equation}
The last two terms in the right hand side cancel each other, leaving a
trivial dependence in space 2 and ${\rm Tr}_1\Big( {L_1^+(q^{\frac
c2}z)}^{t_1} \widetilde{L}_1^-(z)\Big) \equiv t(z)$ in space 1.  This
shows the commutation of $t(z)$ with $L^+(w)$ and therefore with
$L^-(w) = (g^{\frac 12} h g^{\frac 12}) L^+(-p^{\frac {1}{2}}
q^{-\frac{c}{2}}w) (g^{\frac 12} h g^{\frac 12})^{-1}$, hence with the
full algebra ${\cal A}_{q,p}(\widehat{sl}(N)_{c})$ at $c=-N$.
\finproof

\medskip

\noindent This demonstration reproduces the proof for $N=2$ given in 
\cite{AFRS97a}; note that only in the $sl(N)$ crossing-symmetry relation 
(\ref{eq223}) does $N$ appear explicitly.  The form of the operator 
(\ref{eq31}) is identical to the form of the commuting operator derived 
by \cite{RSTS} in the case of the quantum algebra ${\cal 
U}_{q}(\widehat{sl}(N)_{c})$.  As in \cite{AFRS97a} the center of ${\cal 
A}_{q,p}(\widehat{sl}(N)_{c})$ at $c=-N$ may contain other generators 
which we have not yet derived.  However, $t(z)$ do close on their own 
a Poisson algebra as we are going to show.

\subsection{Exchange algebra}

In order to get the Poisson structure on $t(z)$, we need to compute
the exchange algebra between the operators $t(z)$ and $t(w)$ when $c
\ne -N$.  From the definition of the element $t(z)$, one has
\begin{equation}
t(z)t(w) = L(z)^{i_1}_{i_1} \, L(w)^{i_2}_{i_2} =
L^+(q^{\frac c2}z)^{j_1}_{i_1} \, \widetilde L^-(z)^{j_1}_{i_1} \,
L^+(q^{\frac c2}w)^{j_2}_{i_2} \, \widetilde L^-(w)^{j_2}_{i_2} \,.
\label{eq310}
\end{equation}
Suitable rewritings of the relations (\ref{eq236}) lead to the
following exchange relations between the operators $L^+$ and
$\widetilde{L}^-$:
\begin{subequations}
\begin{eqnarray}
&& \widetilde{L}_1^-(z) \, L_2^+(w) = \Big(\widehat R_{21}
(q^{\frac c2}w/z)^{t_1}\Big)^{-1} \, L_2^+(w) \, \widetilde{L}_1^-(z) 
\, \widehat R_{21}^*(q^{-\frac c2}w/z)^{t_1} \,, \\
&& L_1^+(z) \, \widetilde{L}_2^-(w) = \widehat R_{12}(q^{\frac
c2}z/w)^{t_2} \, \widetilde{L}_2^-(w) \, L_1^+(z) \, \Big(\widehat
R_{12}^*(q^{-\frac c2}z/w)^{t_2}\Big)^{-1} \,, \\
&& \widetilde{L}_1^-(z) \, \widetilde{L}_2^-(w) =
\Big(\widehat R_{12}(z/w)^{t_1t_2}\Big)^{-1} \,
\widetilde{L}_2^-(w) \, \widetilde{L}_1^-(z) \, \widehat
R_{12}^*(z/w)^{t_1t_2} \,.
\end{eqnarray}
\label{eq311}
\end{subequations}
The exchange relations (\ref{eq236}), (\ref{eq311}) and
the properties of the matrix $R_{12}$ given in Theorem
\ref{thmone} then allow to move the matrices $L^+(w)$, $\widetilde
L^-(w)$ to the left of the matrices $L^+(z)$, $\widetilde L^-(z)$.  One
obtains
\begin{equation}
t(z)t(w) = {\cal Y}(w/z)^{i_1i_2}_{j_1j_2} ~ L(w)^{j_2}_{i_2} ~
L(z)^{j_1}_{i_1} \,,
\label{eq312}
\end{equation}
where ${\cal Y}(w/z) = T(w/z) {\cal M}(w/z)$ with
\begin{equation}
{\cal M}(w/z) = \left(\left(R_{21}(w/z) \,
{R_{21}(q^{c+N}w/z)}^{-1} \, {R_{12}(z/w)}^{-1}\right)^{t_2}
\, {R_{12}(q^cz/w)}^{t_2}\right)^{t_2}
\label{eq313}
\end{equation}
and
\begin{equation}
T(w/z) = \frac{\tau_N(q^{\frac 12}z/w)\tau_N(q^{\frac 12 - c}w/z)}
{\tau_N(q^{\frac 12}w/z)\tau_N(q^{\frac 12 - c}z/w)} \,.
\label{eq314}
\end{equation}

\subsection{Poisson structures on the center of ${\cal
A}_{q,p}(\widehat{sl}(N)_{c})$ at $c=-N$}

One enounces:
\begin{thm}\label{thmthree}
The elements $t(z)$ form a closed algebra under the natural Poisson
bracket on the center of ${\cal A}_{q,p}(\widehat{sl}(N)_{c})$ given by
(with $x=w/z$)
\begin{equation}
\{ t(z),t(w) \} = - \ln q \left(x\frac{d}{dx}\ln\tau_N(q^{\frac 12}x) - 
x^{-1}\frac{d}{dx^{-1}} \ln\tau_N(q^{\frac 12}x^{-1}) \right) \, t(z) \, 
t(w) \,.
\label{eq323}
\end{equation}
\end{thm}

\medskip

\noindent
{\bf Proof:}
At the critical level $c=-N$, it is easy to show by direct
calculation from (\ref{eq313}) and (\ref{eq314}) that
\begin{equation}
T(x)_{cr} = 1 \quad \mbox{and} \quad {\cal M}(x)_{cr} = \II_2
\otimes \II_2
\label{eq315}
\end{equation}
One recovers immediately that $t(z)t(w) = t(w)t(z)$ at the critical
level $c=-N$.  Hence a natural Poisson structure can be defined by
\begin{equation}
\Big\{ t(z),t(w) \Big\} =
\left(\frac{d{\cal Y}}{dc}(w/z)\right)^{i_1i_2}_{j_1j_2} \,
L(w)^{j_2}_{i_2} \, L(z)^{j_1}_{i_1}\bigg\vert_{cr} \,.
\label{eq316}
\end{equation}
{From} eq. (\ref{eq315}), one has
\begin{equation}
\frac{d{\cal Y}}{dc}(x)\bigg\vert_{cr} =
\frac{dT}{dc}(x)\bigg\vert_{cr} \II_2 \otimes \II_2 \, + \,
\frac{d{\cal M}}{dc}(x)\bigg\vert_{cr} \,.
\label{eq317}
\end{equation}
Now, one has
\begin{equation}
\frac{d{\cal M}}{dc}\bigg\vert_{cr} = -\ln q \left\{
\left(y\frac{d}{dy}R_{21} \bigg\vert_{y=w/z}\right)^{t_2}
R_{12}(q^{-N}z/w)^{t_2} - (R_{12}(z/w)^{-1})^{t_2}
\left(y\frac{d}{dy}R_{12}\bigg
\vert_{y=q^{-N}z/w}\right)^{t_2} \right\}
\label{eq318}
\end{equation}
Taking now the derivative of the relation (\ref{eq230}), one obtains
\begin{equation}
\left( R_{12}(x)^{t_2} \right)^{-1} \left( y\frac{d}{dy}R_{12} 
\bigg\vert_{y=x} \right)^{t_2} \left( R_{12}(x)^{t_2} \right)^{-1} = 
\left( -y\frac{d}{dy}R_{12}^{-1} \bigg\vert_{y=xq^N} \right)^{t_2}
\label{eq319}
\end{equation}
which can be rewritten by virtue of eq. (\ref{eq230}) as:
\begin{equation}
\left( R_{12}(xq^{N})^{-1} \right)^{t_2} \left( y\frac{d}{dy}R_{12} 
\bigg\vert_{y=x} \right)^{t_2} \left( R_{12}(xq^{N})^{-1} \right)^{t_2} 
= \left( -y\frac{d}{dy}R_{12}^{-1} \bigg\vert_{y=xq^N} \right)^{t_2}
\label{eq320}
\end{equation}
Using then eq.  (\ref{eq320}) for the value $x=q^{-N}z/w$ and
inserting it into eq.  (\ref{eq318}), one finds
\begin{equation}
\frac{d{\cal M}}{dc}\bigg\vert_{cr} = -\ln q \left\{ \left( 
y\frac{d}{dy}R_{21} \bigg\vert_{y=w/z} \right)^{t_2} + \left( 
y\frac{d}{dy}R_{12}^{-1} \bigg\vert_{y=z/w} \right)^{t_2} \right\} 
R_{12}(q^{-N}z/w)^{t_2}
\label{eq321}
\end{equation}
{From} the unitarity property, it follows that
$\displaystyle \frac{d{\cal M}}{dc}\bigg\vert_{cr}=0$. Hence
$\displaystyle \frac{d{\cal Y}}{dc}\bigg\vert_{cr} =
\frac{dT}{dc}\bigg\vert_{cr} \II_2 \otimes \II_2$.
This guarantees that the Poisson bracket of $t(z)$ closes on $t(z)$, a
property not obvious since we have indicated that $t(z)$ may not exhaust 
the center of ${\cal A}_{q,p}(\widehat{sl}(N)_{c})$.  \\
Finally, using the $q^{N}$-periodicity property of the function
$\tau_{N}$ and eq. (\ref{eq315}), the derivative of $T(x)$ is given by:
\begin{equation}
\frac{dT}{dc}(x)\bigg\vert_{cr} = \frac{d\ln
T}{dc}(x)\bigg\vert_{cr} = -\ln q
\left(x\frac{d}{dx}\ln\tau_N(q^{\frac 12}x) -
x^{-1}\frac{d}{dx^{-1}} \ln\tau_N(q^{\frac 12}x^{-1}) \right) \,.
\label{eq322}
\end{equation}
\finproof

\subsection{Explicit Poisson structures at $c=-N$}

The structure function in eq.  (\ref{eq323}) is easily computed. 
One obtains:
\begin{equation}
\{ t(z),t(w) \} = f(w/z) \, t(z) \, t(w) \,,
\label{eq324}
\end{equation}
where
\begin{eqnarray}
f(x) &=& -2\ln q \left[ \sum_{\ell \ge 0}
\, \left( \frac{2x^2q^{2N\ell}}{1-x^2q^{2N\ell}} -
\frac{x^2q^{2N\ell+2}}{1-x^2q^{2N\ell+2}} -
\frac{x^2q^{2N\ell-2}}{1-x^2q^{2N\ell-2}} \right ) \right. \nonumber \\
&& \left.  - \frac{x^2}{1-x^2} + \sfrac{1}{2}
\frac{x^2q^2}{1-x^2q^2} + \sfrac{1}{2} \frac{x^2q^{-2}}{1-x^2q^{-2}} 
- (x \leftrightarrow x^{-1}) \right] \,.
\label{eq325}
\end{eqnarray}
We now define the Poisson structure for modes of the generating function 
derived from (\ref{eq324}, \ref{eq325}).  The modes of $t(z)$ are defined 
in the sense of generating functions (or formal series expansions):
\begin{equation}
t_n = \oint_C \frac{dz}{2\pi iz} \, z^{-n} \, t(z)
\label{eq326}
\end{equation}
where $C$ is a contour encircling the origin.  \\
Mode expansions when the structure function $f(x)$ has an infinite set
of poles at $x = q^{P(k)}$ (where $P(k)$ is integer) require a specific
definition using the notion of ``sectors''. This was done at the
classical level in \cite{AFRS97a}.  The procedure runs as follows:  

\medskip

The Poisson bracket between the modes is given by a double contour
integral:
\begin{equation}
\{ t_n,t_m \} = \oint_{C_1} \frac{dz}{2\pi iz} \oint_{C_2} 
\frac{dw}{2\pi iw} z^{-n} \, w^{-m} \, f(w/z) \, t(w) \, t(z) \,.
\label{eq327bis}
\end{equation}
The function $f(x)$ has here simple poles at $x = \pm q^{\pm N\ell}$
and $x = \pm q^{\pm N\ell \pm 1}$. Hence the relative position of the
contours $C_{1}$ and $C_{2}$ must be specified in order to have an
unambiguous result for $(\ref{eq327bis})$.  In addition, the
antisymmetry of the Poisson brackets is only guaranteed at the
mode level by an explicit symmetrization of $(\ref{eq327bis})$ with
respect to the position of the contours $C_{1}$ and $C_{2}$. We shall 
comment on this fact when discussing the quantum mode-exchange structure. 
The mode Poisson bracket is thus defined as:
\begin{equation}
\{ t_n,t_m \} = \frac{1}{2} \left( \oint_{C_1} \frac{dz}{2\pi iz}
\oint_{C_2} \frac{dw}{2\pi iw} + \oint_{C_2}
\frac{dz}{2\pi iz} \oint_{C_1} \frac{dw}{2\pi iw} \right)
z^{-n} \, w^{-m} \, f(w/z) \, t(z) \, t(w) \,,
\label{eq327}
\end{equation}
where $C_{1}$ and $C_{2}$ are circles of radii $R_{1}$ and $R_{2}$ and 
one chooses $R_{1} > R_{2}$.  Explicit evalutation of $(\ref{eq327})$ 
now requires to express $f(w/z)$ as a convergent Laurent series in the 
appropriate domains for $\vert w/z \vert$.  

\medskip

Let us define the sector $(k)$ by $\displaystyle \frac{R_{1}}{R_{2}} \in 
\left] {\vert q\vert }^{-P(k)},{\vert q\vert }^{-P(k+1)} \right[ \,$, 
$P(k)$ being the ordered set of powers of $q^{-1}$ where the poles of 
$f$ are located (here $\{P(k), \, k \in \NN\} = \{0,1; N-1,N,N+1; \dots; 
N\ell-1,N\ell,N\ell+1; \dots\}$ respectively).  For every sector $(k)$, 
eq.  (\ref{eq327}) defines a distinct mode Poisson bracket.

\medskip

As in the $sl(2)$ case, one observes the difference between the
analytic continuation formula (\ref{eq323}), which is unique, and the
formal series formula (\ref{eq327}), where every $k$-labeled
convergent Laurent series expansions for $f$ may be taken as the
formal series expansion of $f$.  This fact is also mentioned in
\cite{FF}, considering the quantum problem.
\begin{prop}
In the case $k=0$ (i.e.  $\displaystyle \frac{R_{1}}{R_{2}} \in 
\left] 1,{\vert q\vert }^{-1} \right[ \,$), one finds
\begin{equation}
\{ t_n,t_m \}_{k=0} = -2\ln q (q-q^{-1}) \sum_{r\in\ZZ}
\frac{[(N-1)r]_{q}[r]_{q}}{[Nr]_{q}} \, t_{n-2r} t_{m+2r} \,,
\label{eq328}
\end{equation}
where the $q$-numbers $[r]_{q}$ are defined as usual:
\begin{equation}
[r]_{q} \equiv \frac{q^r-q^{-r}}{q-q^{-1}} \,.
\label{eq329}
\end{equation}
\end{prop}

\noindent
The proof is immediate.  \\
When $k \ne 0$, one must add to (\ref{eq328}) contributions arising
from the poles at $q^{-P(j)}$ with $j=1,\dots,k$. 
\begin{prop}
The convergent series expansions in any sector $(k)$ are obtained by adding
to the coefficients of the convergent series at $k=0$, coefficients
obtained from the canonical formal series expansion of the distributions
$\delta(q^{-P(j)}w/z)+\delta(-q^{-P(j)}w/z) -\delta(q^{P(j)}w/z)-
\delta(-q^{P(j)}w/z)$ for $j=1,\dots,k$, where $\delta(x) \equiv \sum_{n \in 
\ZZ} x^n$ for $x \in \CC$.
\end{prop}

\medskip

\noindent
{\bf Proof:} Moving from a sector $(k)$ to a sector $(k+1)$ requires to 
rewrite the only term from (\ref{eq325}) whose series expansion becomes 
divergent, namely $\displaystyle 
\frac{x^{2}q^{2P(k)}}{1-x^{2}q^{2P(k)}}$, by a convergent series 
expansion for $\vert x \vert > {\vert q\vert }^{-P(k)}$.  This 
substitutes the series $-\sum_{r \ge 0} x^{-2r}q^{-2P(k)r}$, convergent 
for $\vert x \vert > {\vert q \vert}^{-P(k)}$, to the series $\sum_{r>0} 
x^{2r}q^{2P(k)r}$, convergent when $\vert x \vert <{\vert q 
\vert}^{-P(k)}$.  The overall result, in the full series expansion, is 
to ``add'' the difference (order by order in $x^{2r}$), namely 
$-\sum_{r\in\ZZ} x^{2r}q^{2P(k)r} = -\delta(x^{2}q^{2P(k)})$.  A similar 
reasoning generates the term $\delta(x^{2}q^{-2P(k)})$ from the $x^{-1}$ 
terms in (\ref{eq325}).  
\finproof

\medskip

Beware that the terms obtained at $\vert x\vert ={\vert q \vert}^{\pm 
N \ell}$ get an overall $2$ factor while the terms at $\vert x\vert 
={\vert q\vert}^ {\pm N \ell \pm 1}$ get an overall $-1$ factor.

Specification of a Poisson structure in the context of a 
multiple-singularities structure function therefore does require going 
to an explicit mode expansion.

\subsection{Realization of the higher spin generators}\label{sect35}

To realize deformed ${\cal W}_{N}$ Poisson structures, we need to 
introduce generating functions for the higher spin objects. Having at 
our disposal only one commuting generating function $t(z)$, we are led 
by comparison with \cite{FR} to define shifted products, although with 
the same generator.  Notice that such ordered shifted products were 
used a long time ago, to construct trigonometric and elliptic 
$R$-matrices from rational ones by the so-called ``mean procedure'' 
\cite{Mik79,FadRe83}.  We define accordingly:
\begin{equation}
s_{i}(z) = \prod_{u=-(i-1)/2}^{(i-1)/2} t(q^uz)
\label{eq330}
\end{equation}
where $i=1,\dots,N-1$.  \\
The generators $s_{1}(z),\dots,s_{N-1}(z)$ then close a Poisson algebra
with the following Poisson brackets ($i,j=1,\dots,N-1$):
\begin{equation}
\{ s_{i}(z),s_{j}(w) \} = \sum_{u=-(i-1)/2}^{(i-1)/2}
\sum_{v=-(j-1)/2}^{(j-1)/2}
f\Big(q^{v-u}\frac{w}{z}\Big) \, s_{i}(z) \, s_{j}(w) \,.
\label{eq331}
\end{equation}
Although the generating functions (\ref{eq330}) are here all constructed 
from one single object $t(z)$, we shall see that the Poisson structure 
deduced from (\ref{eq331}) do not reflect this dependence and give rise 
indeed to genuine ${\cal W}_{N}$-type structures, in particular 
recovering the $q$-${\cal W}_{N}$ algebra in \cite{FR}, as a 
consequence of the sector structure of the $q$-${\cal W}_{N}$ algebra 
in terms of modes. Notice also from (\ref{eq331}) that the limitation 
of the index $i$ to values smaller than $N$ is justifed by the fact that 
a product of $N$ functions $t(q^uz)$ in (\ref{eq330}) in fact Poisson 
commute with all $s_{i}(w)$ owing to:
\begin{equation}
\sum_{u=1}^{N} f(q^ux) = 0 \,, \qquad \forall \, x \in \CC
\label{eq331bis}
\end{equation}
We now study the mode expansion of (\ref{eq331}).  The singularities of 
the structure function lie at $x =\pm q^{\pm N\ell}q^u$ and $\pm 
q^{\pm N\ell \pm 1}q^u$ with $\ell$ a positive integer and $u$ an 
integer (resp.  a half-integer) from $1-\sfrac{1}{2}(i+j)$ to 
$\sfrac{1}{2}(i+j)-1$ for $(i+j)$ even (resp.  $(i+j)$ odd).  They fall 
into sets of $(i+j+1)$ poles symmetrically arranged around $q^{\pm 
N\ell}$, separated by one power of $q$.  This setting defines a labeling 
of sectors for the Poisson brackets as follows: for fixed $i$ and $j$, 
we define a sector $(k)$ by $\displaystyle \frac{R_{1}}{R_{2}} \in 
\left] {\vert q\vert }^{-P_{ij}(k)},{\vert q\vert }^{-P_{ij}(k+1)} 
\right[ \,$, $P_{ij}(k)$ being the ordered set of powers of $q^{-1}$, 
positive and negative, where the poles of the structure function of 
(\ref{eq331}) are located, such that $P_{ij}(0) = 0$.  
The rules which define the sectors in which the Poisson bracket of the 
modes of a given couple $(s_{i},s_{j})$ is computed, are the following:  \\
$\bullet$ Poisson brackets between modes of the same field
$\{s_{i}(z),s_{i}(w)\}$ are required by antisymmetry to be computed
on contours $C_{1}$ and $C_{2}$ symmetrized as in formula (\ref{eq327}).
Hence they are labeled by positive numbers $(k)$ only, corresponding to the
choice $\displaystyle \frac{R_{1}}{R_{2}} \in \left]
{\vert q\vert}^{-P_{ii}(k)},{\vert q\vert}^{-P_{ii}(k+1)} \right[$.  \\
$\bullet$ Poisson brackets between modes of different fields 
$\{s_{i}(z),s_{j}(w)\}$ can be computed on a single set of contours 
$C_{1}$ and $C_{2}$ such that $\displaystyle \frac{R_{1}}{R_{2}} \in 
\left] {\vert q\vert}^{P_{ij}(k)},{\vert q\vert}^{P_{ij}(k+1)} \right[$ where 
$P_{ij}(k)$ may be positive or negative.  Symmetrization over $C_{1}$ and 
$C_{2}$ is not required.  Antisymmetry of the Poisson bracket is imposed 
by computing $\{s_{i},s_{j}\}_{(k)}$ for $i<j$ and {\em setting} 
$\{s_{j},s_{i}\}_{(k)} \equiv -\{s_{i},s_{j}\}_{(k)}$.  Hence these Poisson 
brackets are labeled by positive and negative numbers, one for each 
couple $(i,j)$.  \\
$\bullet$ The choice of sectors $k(i,j)$ on which Poisson brackets of 
different couples are computed is arbitrary.  In fact, quadratic Poisson 
bracket structures obey the Jacobi identity as soon as they are 
antisymmetric, hence any antisymmetric Poisson structure is consistent.  \\
To summarize, a complete Poisson structure for $\{s_{i},s_{j}\}$ is
characterized by the choice of $N-1$ positive integer labels and
$\sfrac{1}{2} (N-1)(N-2)$ integer labels.
\begin{prop}
When one chooses $k=0$ for all sectors, one obtains a compact generic 
expression (with $n,m \in \ZZ$):
\begin{equation}
\{ s_i(n),s_j(m) \} = -2\ln q (q-q^{-1}) \sum_{r\in\ZZ}
\frac{[(N-\max(i,j))r]_{q}[\min(i,j)r]_{q}}{[Nr]_{q}}
\, s_i(n-2r) s_j(m+2r) \,.
\label{eq332}
\end{equation}
\end{prop}
The proof of this formula is given in Appendix B. This Poisson bracket 
structure is identical to the Poisson bracket structure obtained in 
\cite{FR} from a bosonization construction, excluding the extra 
$\delta$-type terms in $s_{i-p}s_{j+p}$.  We shall comment on the 
possibility of occurence for such terms in the conclusion.  This 
realizes, as in the $sl(2)$ case, a non-trivial connection between the 
$q$-${\cal W}_{N}$ algebra and the $sl(N)$ elliptic algebra.
\begin{prop}
Any Poisson structure in a given sector $k(i,j)$ can be obtained from 
(\ref{eq332}) by adding to the $r$-dependent structure coefficient 
contributions from the relevant singularities of the structure function.  
They are given by formal power series expansions of terms $\delta(\pm 
q^{\pm P_{ij}(s)}w/z) s_{i}(z) s_{j}(w)$ for $s=1,\dots,k(i,j)$.
\end{prop}

\medskip

One should now emphasize that the Poisson structures (\ref{eq332}) in 
any sector $k(i,j)$ are not identified to the structure which would be 
obtained from application by Leibniz rule to the mode expansion 
$\sum_{a_{1}+\dots+a_{i}=m} \prod t_{a_{i}}q^{k_{i}}$ obtained from 
(\ref{eq330}) by a single contour integral on a contour $C_{1}$ for $z$ 
of any particular $k$-sector Poisson structure for the generators $t_m$ 
derived in (\ref{eq328}).  Indeed, this structure would simply be given 
by the corresponding structure function $\displaystyle 
\frac{[(N-1)r]_{q}[r]_{q}}{[Nr]_{q}}$ plus its $\delta$-contributions.

By contrast, it can be seen in (\ref{eq331})-(\ref{eq332}) and in the 
derivation (Appendix B) that (forgetting for the time being the further 
symmetrization requirements over the double contour integral) the 
Poisson structure (\ref{eq332}) follows in fact from application of the 
Leibniz rule to expansions of the form above, but where each individual 
Poisson bracket $\{t_{n},t_{m}\}$ must be computed in {\em distinct} 
relative sectors since they stem from Poisson brackets between 
generating functions $t(q^{u_1}z)$ and $t(q^{u_2}w)$ given by 
$f(q^{u_1-u_2}z/w)$ as a contribution to the structure function, where 
$u_{1}$ and $u_{2}$ respectively live in two intervals 
$\left]-\sfrac{1}{2}(i-1),\sfrac{1}{2}(i-1)\right[$ and 
$\left]-\sfrac{1}{2}(j-1),\sfrac{1}{2}(j-1)\right[$ : 
hence the relative position of integration contours depend on 
the difference between the indices $u_{1}$ and $u_{2}$ which lives in 
$\left]-\sfrac{1}{2}(i+j)-1,\sfrac{1}{2}(i+j)-1\right[$
and no individual Poisson bracket $\{t_{n},t_{m}\}$ can be 
factored out.  This is in particular true when considering identical 
fields $\{s_{i},s_{i}\}$.  The symmetrization procedure required for 
the $t_{n}$ Poisson brackets adds a further obstacle to attempts at 
factoring out {\em symmetrized} Poisson brackets for the modes $t_{n}$.

To summarize, once the Poisson bracket of composite fields 
$s_{i}(z),s_{j}(w)$ are computed for the modes defined by contour 
integrals in specified relative positions for $z$ and $w$, giving 
(\ref{eq331})-(\ref{eq332}), the nature of composite fields $s_{i}(z)$ 
as products of the initial $t(z)$ generators (\ref{eq330}) is 
obliterated from the new mode Poisson structure thus obtained.  
The composite fields then assume the nature of independent objects with 
the Poisson structure (\ref{eq332}), thereby validating the seemingly 
redundant definition (\ref{eq330}).

\medskip

{\bf Remark:} Reciprocally, all supplementary terms denoted $(\delta(\pm 
q^uw/z) s_{i}(z) s_{j}(w))_{n,m}$ from the mode expansion in $k \ne 0$ 
sectors are defined in the sense of {\em formal series} expansions as:
\[
(\delta(\pm q^uw/z) s_{i}(z) s_{j}(w))_{n,m} \equiv
\sum_{p \in \ZZ} q^{up} s_{i}(n-p) s_{j}(m+p) \,.
\]
Their Poisson brackets must therefore be computed consistently by 
Leibniz rule applied to their mode expansion, and {\em not} by using 
$s_{i}(z) s_{j}(q^{u}z)$ as a generating functional.  In particular, the 
extra terms $\delta(q^{-1}z/w) t(z) t(w)$ in the $sl(2)$ case must not 
be understood as a central extension although the generating function 
$t(z) t(zq)$ Poisson commutes with $t(w)$.  Central and lower-spin terms 
do not occur in our derivation.  We shall comment on their absence 
here, and their possible reconstruction, in the conclusion.

\section{Quadratic algebras in ${\cal A}_{q,p}(\widehat{sl}(N)_{c})$}
\label{sect4}
\setcounter{equation}{0}

\begin{thm}\label{thmfour}
In the three-dimensional parameter space generated by $p,q,c$, one
defines a two-dimen\-sional surface $\Sigma_{N,M}$ for any integer $M
\in \ZZ$ by the set of triplets $(p,q,c)$ connected by the relation
$(-p^{\frac 12})^{NM} = q^{-c-N}$.  On the surface $\Sigma_{N,M}$, the
generators $t(z)$ realize an exchange algebra with the generators
$L(w)$ of ${\cal A}_{q,p}(\widehat{sl}(N)_{c})$:
\begin{equation}
t(z) \, L(w) = F\Big(M,\frac{w}{z}\Big) \, L(w) \, t(z)
\label{eq41}
\end{equation}
where
\begin{subequations}
\begin{eqnarray}
F(M,x) &=& q^{2M(N-1)} \prod_{k=0}^{NM-1}
\frac{\Theta_{q^{2N}}(x^{-2}p^{-k}) \, \Theta_{q^{2N}}(x^2p^{k})}
{\Theta_{q^{2N}}(x^{-2}q^2p^{-k}) \, \Theta_{q^{2N}}(x^2q^2p^{k})}
\quad \mbox{for $M>0$} \,, \\
F(M,x) &=& q^{-2|M|(N-1)} \prod_{k=1}^{N|M|}
\frac{\Theta_{q^{2N}}(x^{-2}q^2p^{k}) \, \Theta_{q^{2N}}(x^2q^2p^{-k})}
{\Theta_{q^{2N}}(x^{-2}p^{k}) \, \Theta_{q^{2N}}(x^2p^{-k})}
\quad \mbox{for $M<0$} \,.
\end{eqnarray}
\label{eq42}
\end{subequations}
\end{thm}

\medskip

\noindent
{\bf Proof:} The proof runs along similar lines to the commutativity
proof of Theorem \ref{thmtwo}.  From eqs.  (\ref{eq31}) to (\ref{eq34}),
one gets
\begin{equation}
t(z) \, L_2^+(w) = {\rm Tr}_1\Big(L_1^+(q^{\frac c2}z)^{t_1}
(\widehat R_{21}(q^{\frac c2}w/z)^{t_1})^{-1} L_2^+(w) \widetilde
L_1^-(z) \widehat R_{21}^*(q^{-\frac c2}w/z)^{t_1}\Big) \,.
\label{eq43}
\end{equation}
One also has, from (\ref{eq36}):
\begin{equation}
L_1^+(q^{\frac c2}z)^{t_1} \, (\widehat R_{21}(q^{\frac c2 
-c-N}w/z)^{t_1})^{-1} \, L_2^+(w) = L_2^+(w) \, (\widehat 
R_{21}^*(q^{-\frac c2}w/z)^{-1})^{t_1} \, L_1^+(q^{\frac c2}z)^{t_1} \,,
\label{eq44}
\end{equation}
One realizes that the only obvious condition that allow a substitution
of eq.  (\ref{eq44}) into eq.  (\ref{eq43}) using the
quasi-periodicity of the matrix $\widehat R_{12}$ is the following:
\begin{equation}
(-p^{\frac 12})^{NM} = q^{-c-N} \quad \mbox{with} \,\, M \in \ZZ
\label{eq45}
\end{equation}
Actually, from the quasi-periodicity property of $\widehat R_{12}$, one
has:
\begin{equation}
\widehat R_{21}((-p^{\frac 12})^{NM}x) = F(M,x) \widehat R_{21}(x)
\label{eq46}
\end{equation}
where
\begin{equation}
F(M,x) = \left\{ \begin{array}{ll}
\displaystyle \prod_{k=0}^{NM-1} F\left((-p^{\frac 12})^kx\right)
& \mbox{for $M>0$} \,, \\ \\
\displaystyle \prod_{k=1}^{N|M|} F\left((-p^{\frac 12})^{-k}x
\right)^{-1}
& \mbox{for $M<0$} \,, \\
\end{array} \right.
\label{eq47}
\end{equation}
the function $F(x)$ being given by
\begin{equation}
F(x) = \tau_{N}^{-1}(q^{\frac 12}x)
\tau_{N}^{-1}(q^{\frac 12} x^{-1}) = q^{2-2/N}
\frac{\Theta_{q^{2N}}(x^2) \Theta_{q^{2N}}(x^{-2})}
{\Theta_{q^{2N}}(q^2x^2) \Theta_{q^{2N}}(q^2x^{-2})} \,.
\label{eq48}
\end{equation}
Then on the two-dimensional surface defined by (\ref{eq45}), the
equation (\ref{eq44}) becomes
\begin{equation}
L_1^+(q^{\frac c2}z)^{t_1} \, (\widehat R_{21}(q^{\frac 
c2}w/z)^{t_1})^{-1} \, L_2^+(w) = F\Big(M,q^{\frac c2}\frac{w}{z}\Big) 
L_2^+(w) \, (\widehat R_{21}^*(q^{-\frac c2}w/z)^{-1})^{t_1} \, 
L_1^+(q^{\frac c2}z)^{t_1} \,.
\label{eq49}
\end{equation}
It follows that
\begin{equation}
t(z) \, L_2^+(w) = F\Big(M,q^{\frac c2}\frac{w}{z}\Big) L_2^+(w)
{\rm Tr}_1\Big((\widehat R_{21}^*(q^{-\frac c2}w/z)^{-1})^{t_1} \,
L_1^+(q^{\frac c2}z)^{t_1} ({L_1^-(z)}^{-1})^{t_1}
\widehat R_{21}^*(q^{-\frac c2}w/z)^{t_1}\Big) \,.
\label{eq410}
\end{equation}
The two $R$ matrices cancel due to the relation ${\rm Tr}_1 \Big(
R_{21} Q_1 {R'}_{21} \Big) = {\rm Tr}_1 \Big( Q_1 {R'_{21}}^{t_2}
{R_{21}}^{t_2} \Big)^{t_2}$.  Hence recalling that $L^+(w) =
L(q^{\frac c2}w)$, one gets the desired result.  Finally, one needs to
compute the factor $F(M,x)$.
Using eq.  (\ref{eq48}), one obtains the expressions of Theorem
\ref{thmfour}.
\finproof

\medskip

\noindent
\begin{thm}\label{thmfive}
On the surface $\Sigma_{N,M}$, $t(z)$ closes a quadratic subalgebra:
\begin{equation}
t(z)t(w) = {\cal Y}_{N,p,q,M}\Big(\frac{w}{z}\Big) \, t(w)t(z)\, ,
\label{eq411}
\end{equation}
where
\begin{equation}
{\cal Y}_{N,p,q,M}(x) = \left\{ \begin{array}{ll}
\displaystyle \prod_{k=1}^{NM} \frac{\Theta_{q^{2N}}^2(x^2 p^{-k}) \,
\Theta_{q^{2N}}(x^2 q^2 p^k) \, \Theta_{q^{2N}}(x^2 q^{-2} p^k)}
{\Theta_{q^{2N}}^2(x^2 p^k) \, \Theta_{q^{2N}}(x^2 q^2 p^{-k})
\, \Theta_{q^{2N}}(x^2 q^{-2} p^{-k})}
& \mbox{for $M>0$} \,, \\ \\
\displaystyle \prod_{k=1}^{N|M|-1}
\frac{\Theta_{q^{2N}}^2(x^2 p^{-k}) \,
\Theta_{q^{2N}}(x^2 q^2 p^k) \, \Theta_{q^{2N}}(x^2 q^{-2} p^k)}
{\Theta_{q^{2N}}^2(x^2 p^k) \, \Theta_{q^{2N}}(x^2 q^2 p^{-k})
\, \Theta_{q^{2N}}(x^2 q^{-2} p^{-k})}
& \mbox{for $M<0$} \,. \\
\end{array} \right.
\label{eq412}
\end{equation}
\end{thm}

\medskip

\noindent
{\bf Proof:} From Theorem \ref{thmfour}, one has
\begin{subequations}
\begin{eqnarray}
t(z) \, L^+(w) &=& F\Big(M,q^{\frac c2}\frac{w}{z}\Big)
\, L^+(w) \, t(z) \,, \\
t(z) \, (L^-(w))^{-1} &=& F^{-1}\Big(M,-p^{\frac 12}\frac{w}{z}\Big)
\, (L^-(w))^{-1} \, t(z) \,,
\end{eqnarray}
\label{eq413}
\end{subequations}
Hence, the definition (\ref{eq31}) of $t(z)$ immediately implies:
\begin{equation}
t(z) \, t(w)= \frac{\displaystyle F\Big(M,q^{c}\frac{w}{z}\Big)}
{\displaystyle F\Big(M,-p^{\frac 12}\frac{w}{z}\Big)} \,\,
t(w) \, t(z) \,.
\label{eq414}
\end{equation}
The explicit expression (\ref{eq42}) for $F(M,x)$ gives the result
as stated above.
\finproof

\medskip

\noindent
We shall discuss the mode expansion of (\ref{eq411}) in Sect.
\ref{sect6}.

\section{Poisson structures associated to commuting subalgebras in
${\cal A}_{q,p}(\widehat{sl}(N)_{c})$}
\label{sect5}
\setcounter{equation}{0}

\begin{thm}\label{thmsix}
On the surface $\Sigma(N,M)$, when $p = q^{Nh}$ with $h \in
\ZZ \backslash \{ 0 \}$, the function ${\cal Y}_{N,p,q,M}$ is equal
to 1. Hence $t(z)$ realizes an Abelian subalgebra in ${\cal
A}_{q,p}(\widehat{sl}(N)_{c})$.
\end{thm}

\medskip

\noindent
{\bf Proof:} Theorem \ref{thmsix} is easily proved using the explicit
expression for $F(M,x)$ and the periodicity properties of the
$\Theta_{q^{2N}}$ functions.
\finproof

\medskip

\noindent
{\bf Remark:} Except in the case $N=2$ (see \cite{AFRS97b}), no value of 
$h$ allows here for $t(z)$ to be an element of a possibly extended 
center of ${\cal A}_{q,p}(\widehat{sl}(N)_{c})$.

\medskip

The result of Theorem \ref{thmsix} nevertheless allows us to define 
Poisson structures on the corresponding Abelian algebras.  They are 
obtained as limits of the exchange algebra (\ref{eq411}) when $p = 
q^{Nh}$ with $h \in \ZZ \backslash \{ 0 \}$.  \\
Conversely it follows that (\ref{eq414}) realizes a natural quantization 
of the Poisson structures obtained by this limit, since it realizes an 
``intermediate'' closed exchange algebra, contrary to the situation at 
$c=-N$ where $t(z)$ immediately lies in the center.  This stands in 
contrast with the construction in \cite{FR} where the quantized 
$q$-${\cal W}_{N}$ algebras must be reconstructed by an independent 
quantization of the deformed classical bosons in the Cartan algebra.  We 
see here the key role of the initial 3-parameter structure ${\cal 
A}_{q,p}(\widehat{sl}(N)_{c})$ compared to the 2-parameter quantum 
algebra ${\cal U}_q(\widehat{sl}(N)_{c})$ used in \cite{FR}.  It allows 
for an intermediate 2-parameter step at $(-p^\sfrac{1}{2})^{NM} = 
q^{-c-N}$ where the generators $t(z)$ themselves close an exchange 
algebra.  Hence it provides at the same time the classical $q$-deformed 
${\cal W}_{N}$ algebra and its $(q,p)$-deformed quantization. 
\begin{thm}\label{thmseven}
Setting $q^{Nh} = p^{1-\beta}$ for any integer $h \ne 0$, the
$h$-labeled Poisson structure defined by:
\begin{equation}
\{ t(z) , t(w) \}^{(h)} = \lim_{\beta \rightarrow 0} \frac{1}{\beta}
\, \Big(t(z)t(w) - t(w)t(z)\Big)
\label{eq51}
\end{equation}
has the following expression:
\begin{equation}
\{ t(z) , t(w) \}^{(h)} = f_h(w/z) \, t(z) \, t(w)
\label{eq52}
\end{equation}
where
\begin{subequations}
\begin{eqnarray}
f_h(x) &=&
2Nh \ln q \left[ \sum_{\ell \ge 0}
E(\sfrac{NM}{2})(E(\sfrac{NM}{2})+1)
\left( \frac{2x^2q^{2N\ell}}{1-x^2q^{2N\ell}}
- \frac{x^2q^{2N\ell+2}}{1-x^2q^{2N\ell+2}}
- \frac{x^2q^{2N\ell-2}}{1-x^2q^{2N\ell-2}} \right) \right.
\nonumber \\
&& + E(\sfrac{NM+1}{2})^2 \left(
\frac{2x^2q^{2N\ell+N}}{1-x^2q^{2N\ell+N}}
- \frac{x^2q^{2N\ell+N+2}}{1-x^2q^{2N\ell+N+2}}
- \frac{x^2q^{2N\ell+N-2}}{1-x^2q^{2N\ell+N-2}} \right)
\nonumber \\
&& \left. - \sfrac{1}{2} E(\sfrac{NM}{2})(E(\sfrac{NM}{2})+1)
\left( \frac{2x^2}{1-x^2} - \frac{x^2q^2}{1-x^2q^2}
- \frac{x^2q^{-2}}{1-x^2q^{-2}} \right)
- (x \leftrightarrow x^{-1}) \right]
\qquad \mbox{for $h$ odd} \,, \nonumber \\
&& \label{eq53a} \\
&& \nonumber \\
&=& N^2M(NM+1)h \ln q \left[ \sum_{\ell \ge 0} \left(
\frac{2x^2q^{2N\ell}}{1-x^2q^{2N\ell}}
- \frac{x^2q^{2N\ell+2}}{1-x^2q^{2N\ell+2}}
- \frac{x^2q^{2N\ell-2}}{1-x^2q^{2N\ell-2}} \right) \right.
\nonumber \\
&& \left. - \sfrac{1}{2} \left( \frac{2x^2}{1-x^2}
- \frac{x^2q^2}{1-x^2q^2} - \frac{x^2q^{-2}}{1-x^2q^{-2}} \right)
- (x \leftrightarrow x^{-1}) \right]
\qquad \mbox{for $h$ even} \,. \nonumber \\
&& \label{eq53b}
\end{eqnarray}
\label{eq53}
\end{subequations}
Here the notation $E(n)$ means the integer part of the number $n$.
\end{thm}
{\bf Proof:} direct calculation.
\finproof

\medskip

\noindent 
The factors $2Nh$ for $h$ odd and $N^2Mh(NM+1)$ for $h$ even are 
inessential and can be reabsorbed into the definition of the classical 
limit as $\beta \rightarrow -Nh\beta$ for $h$ odd and $\beta \rightarrow 
-\sfrac{1}{2} N^2M(NM+1)h\beta$ for $h$ even.  Provided this 
redefinition is done, the formula (\ref{eq53b}) in the case $h$ even 
coincides exactly with the Poisson structure (\ref{eq325}) of the center 
of ${\cal A}_{q,p}(\widehat {sl}(N)_{-N})$.  The formula (\ref{eq53a}) 
in the contrary in the case $h$ odd gives rise to a new Poisson structure.  \\
As before, the equations (\ref{eq53}) lead to infinite families of
Poisson brackets for the modes of $t(z)$ defined by (\ref{eq326}).
For $h$ even, the labeling is identical to the $c=-N$ case. For $h$
odd, the singularities lie at $\pm q^{\pm N\ell}$, $\pm q^{\pm N\ell \pm 1}$
and $\pm q^{\pm N\ell \pm \sfrac{1}{2} N}$, 
$\pm q^{\pm N\ell \pm \sfrac{1}{2} N \pm 1}$, giving rise to extra
triplets of poles lying halfway between the initial ones.  
\begin{prop}
One gets for the Poisson bracket in the $k=0$ sector:
\begin{subequations}
\begin{eqnarray}
\{ t_n , t_m \}^{(h)} &=&
-2\ln q (q-q^{-1}) \sum_{r \in \ZZ} \left(
E(\sfrac{NM}{2})(E(\sfrac{NM}{2})+1)
\frac{[(N-1)r]_{q}[r]_{q}}{[Nr]_{q}} \right. \nonumber \\
&& \hspace{40mm} \left. - E(\sfrac{NM+1}{2})^2 \frac{[r]_{q}^2}{[Nr]_{q}} 
\right) t_{n-2r} t_{m+2r} \hspace{10mm} \mbox{for $h$ odd} \,,
\label{54a} \\
&& \nonumber \\
&=& -2\ln q (q-q^{-1}) \sum_{r \in \ZZ}
\frac{[(N-1)r]_{q}[r]_{q}}{[Nr]_{q}}
t_{n-2r} t_{m+2r} \hspace{24mm} \mbox{for $h$ even} \,.
\label{eq54b}
\end{eqnarray}
\label{eq54}
\end{subequations}
\end{prop}

A realization of higher spin generators is again achieved by 
the formula (\ref{eq330}) with $i=1,\dots,N-1$.  Its justification will 
be the same as in Sect. \ref{sect35}.  The generators
$s_{1}(z),\dots,s_{N-1}(z)$ close a Poisson algebra with the following
Poisson brackets ($i,j=1,\dots,N-1$):
\begin{equation}
\{ s_{i}(z),s_{j}(w) \}^{(h)} = \sum_{u=-(i-1)/2}^{(i-1)/2}
\sum_{v=-(j-1)/2}^{(j-1)/2}
f_h\Big(q^{v-u}\frac{w}{z}\Big) \, s_{i}(z) \, s_{j}(w) \,.
\label{eq55}
\end{equation}
The singularity structure here is as follows.  Singularities of the
function $f_h(x)$ occur at $x =\pm q^{\pm N\ell+u}$,
$\pm q^{\pm N\ell \pm 1+u}$ where $\ell \in \NN$ and $u \in \left[
1-\sfrac{1}{2}(i+j),\sfrac{1}{2}(i+j)-1 \right]$ when $h$ is even.
Additional singularities occur halfway between those ones, at $\pm q^{\pm
N\ell\pm\sfrac{1}{2}N+u}$, $\pm q^{\pm N\ell\pm\sfrac{1}{2}N \pm 1+u}$ when
$h$ is odd.  The sector structure for the Poisson brackets of
$\{s_{i}(z),s_{j}(w)\}$ is easily deduced from these results.  

Easiest to compute are the Poisson brackets in the sector $k=0$ for all 
couples of indices $(i,j)$.  A simpler Poisson structure in this sector 
is defined by furthermore taking in {\em all} Poisson brackets a {\em 
symmetrized} double contour integral with $\displaystyle 
\frac{R_{1}}{R_{2}} \in \left] 1,\vert q \vert^{-1/2} \right[$ (or 
$\left] 1,\vert q \vert^{-1} \right[$ depending on the parity of $N$ and 
$i+j$).  (Note that the symmetrized form is actually not required here 
either when $i \ne j$, but it leads to nicer formulae).  One gets:
\begin{subequations}
\begin{eqnarray}
\{ s_i(n),s_j(m) \}^{(h)} &=& -2\ln q (q-q^{-1}) \sum_{r\in\ZZ} 
\frac{[(N-\max(i,j))r]_{q}[\min(i,j)r]_{q}}{[Nr]_{q}} 
s_i(n-2r) s_j(m+2r) \nonumber \\
&& \hspace{70mm} \mbox{for $h$ even} \,, \label{eq56a} \\
&& \nonumber \\
\{ s_i(n),s_j(m) \}^{(h)} &=& -2\ln q (q-q^{-1}) \sum_{r\in\ZZ} 
\left( E(\sfrac{NM}{2})(E(\sfrac{NM}{2})+1)
\frac{[(N-\max(i,j))r]_{q}[\min(i,j)r]_{q}}{[Nr]_{q}} \right. \nonumber \\
&& \hspace{40mm} \left. - E(\sfrac{NM+1}{2})^2 
\frac{[ir]_{q}[jr]_{q}}{[Nr]_{q}} \right) s_i(n-2r) s_j(m+2r) \nonumber \\
&& \hspace{70mm} \mbox{for $h$ odd and $i+j \le N$} \,, \label{eq56b} \\
&& \nonumber \\
\{ s_i(n),s_j(m) \}^{(h)} &=& -2\ln q (q-q^{-1}) \sum_{r\in\ZZ} 
\left( E(\sfrac{NM}{2})(E(\sfrac{NM}{2})+1)
\frac{[(N-\max(i,j))r]_{q}[\min(i,j)r]_{q}}{[Nr]_{q}} \right. \nonumber \\
&& \left. - E(\sfrac{NM+1}{2})^2 \frac{[ir]_{q}[jr]_{q}}{[Nr]_{q}}
- E(\sfrac{NM+1}{2})^2 [(N-i-j)r]_{q} \right) s_i(n-2r) s_j(m+2r) 
\nonumber \\
&& \hspace{70mm} \mbox{for $h$ odd and $i+j > N$} \,. \label{eq56c}
\end{eqnarray}
\label{eq56}
\end{subequations}
The Poisson bracket (\ref{eq56a}) is again identical to the ``core'' 
contribution in \cite{FR} (i.e.  without lower-spin extensions).  
(\ref{eq56b}-c) however is a completely new type of quadratic 
$q$-deformed classical ${\cal W}_{N}$-algebra with the two extra terms.  
It would be interesting to know about its possible explicit 
constructions.

\section{Quantum exchange algebra}
\label{sect6}
\setcounter{equation}{0}

The exchange algebras (\ref {eq411}) are now understood as natural 
quantizations of the classical $q$-deformed ${\cal W}_N$ algebras, 
including the initial algebra \cite {FR}.  Moreover it follows that the 
exchange algebra (\ref{eq411}) now provides us with building blocks for 
new deformed quantum ${\cal W}_{q,p}(sl(N))$ algebras with an extra 
integer parameter $M$.

We wish to describe here an explicit formulation in terms of modes of 
quantum generating operators $s_i(z)$ defined as in (\ref {eq330}), 
although with a required notion of ordering between individual 
$t(z)$-generators:
\begin{equation}
s_i(z) = \prod_{\frac{i-1}{2} \ge u \ge -\frac{i-1}{2}}^{{\textstyle 
\curvearrowleft}}t(q^uz) \,.
\label{eq61}
\end{equation}
Again justification of this definition as giving genuine $q$-${\cal W}_N$
algebras will come from the arising of sectors in the individual
$t(z)$-$t(w)$ exchange algebra, which eventually combine in a
non-trivial way in the product formula to give rise to a new algebraic
structure.

The exchange algebra from (\ref{eq61}) and (\ref{eq411}) takes the
form:
\begin{equation}
s_i(z) s_j(w) = \prod_{u=-\frac{i-1}{2}}^{\frac{i-1}{2}}
\prod_{v=-\frac{j-1}{2}}^{\frac{j-1}{2}}{\cal Y}_{N,p,q,M} 
\left(q^{v-u}\frac{w}{z}\right) s_j(w) s_i(z) \, ,
\label{eq62}
\end{equation}
for any choice of ordering in (\ref {eq61}). Furthermore it follows
from (\ref {eq61}) and (\ref {eq412}) that $s_N(z)$ commutes with all
other generators; hence we are justified in restricting $i$ to $1,
\ldots , N-1$.

\medskip

Once we choose an exchange function in (\ref {eq412}) by choosing the
integer $M$, the next step in the procedure consists in factorizing
the exchange functions in (\ref {eq62}) into a function analytic
around $w/z=0$ and a function analytic around $z/w=0$. More precisely,
and following \cite {FR2}, one defines a Riemann problem:
\begin{equation}
{\cal Y}_{N,p,q,M;(i,j)}(x)={\cal Y}_+(x){\cal Y}_-^{-1}(x^{-1}) \,
\label{eq63}
\end{equation}
in the neighborhood of a circle $C$ of radius $R$. ${\cal Y}_+$ and
${\cal Y}_-$ are respectively analytic for $|x|<R$ and $|x|>R$.

Varying the values of $R$ with respect to the position of zeroes and 
poles of ${\cal Y}$ leads to different factorizations; moreover even 
with fixed $R$ the solution of (\ref{eq63}) is not unique.  This leads 
to a large choice of acceptable factorizations.  It may be possible to 
choose ${\cal Y}_+={\cal Y}_-$ in the sense of analytic continuation, 
for instance for $i=j=1$, when $R$ takes any value between $|p|^{\frac 
{1}{2}}$ and $|p|^{-\frac{1}{2}}$; this case is treated in detail below.

\medskip

The third step is specific to our approach. It consists in promoting
the exchange relation deduced from (\ref {eq62}) and (\ref {eq63}) to
the level of an analytic extension to the full complex plane of $z/w$:
\begin{equation}
{\cal Y}_- \left(\frac{z}{w}\right) t(z) t(w) = {\cal Y}_+ \left( 
\frac{w}{z}\right) t(w) t(z) \, .
\label{eq64}
\end{equation}
Singularities of ${\cal Y}_\pm$ now play a crucial role as in the
classical case, when we define mode expansions of (\ref {eq64}). The
fields $s_i(z)$ are considered as abstract generating operator-valued
functionals with modes defined by contour integrals:
\begin{equation}
s_i(n)=\oint _C\frac {dz}{2\pi i z}z^{-n}s_i(z) \, .
\label{eq65}
\end{equation}
They obey no particular supplementary relations, which may eventually
arise from explicit realizations of $t(z)$ as in \cite {FF,FR2}.

As a consequence, we need to introduce a similar notion of
``sectors'', determined by the singularities of ${\cal Y}_\pm$, and
defined as the regions in the complex plane for $z/w$ where ${\cal
Y}_+$ and ${\cal Y}_-$ are given by a particular convergent series
expansion (necessarily unique in {\it each} domain since ${\cal Y}_\pm$
are meromorphic).  A choice of relative positions for the contours 
$C_{1}$ for $z$ and $C_{2}$ for $w$ thus gives the unique formal series 
expansions for ${\cal Y}_+$ and ${\cal Y}_-$ to be inserted in the 
double contour integral, which eventually gives an exchange relation for 
the modes.  We first illustrate this on the simplest example $i=j=1$ (i.e. 
$s_{i}(z) = s_{j}(z) \equiv t(z)$).  

\medskip

To extract from (\ref {eq411}) an exchange relation between the
modes $t_n$ (\ref{eq326}) of $t(z)$, we rewrite (\ref{eq411}) in the form:
\begin{equation}
f_{N,p,q,M} \left( \frac{z}{w} \right) \, t(z) \, t(w) =
t(w) \, t(z) \, f_{N,p,q,M} \left( \frac{w}{z} \right) \,.
\label{eq66}
\end{equation}
A consistent choice for $f_{N,p,q,M}(x)$ is then:
\begin{equation}
f_{N,p,q,M}(x) = \prod_{k=1}^{S(M)}(1-x^2)^2 \prod _{n=0}^{\infty} \frac
{(1-x^2p^{-k}q^{2N(1+n)})^2(1-x^2p^kq^{2+2Nn})(1-x^2p^kq^{-2+2N(1+n)})}
{(1-x^2p^kq^{2Nn})^2(1-x^2p^{-k}q^{-2+2N(1+n)})(1-x^2p^{-k}q^{2+2Nn})} \,,
\label{eq67}
\end{equation}
where $S(M)=NM$ for $M>0$ and $S(M)=N|M|-1$ for $M<0$.  \\
In the hypothesis:
\begin{equation}
|p|<1 \,, \quad |q|<1 \,, \quad |p|^{S(M)+1}>|q|^2 \,,
\label{eq68}
\end{equation}
$f_{N,p,q,M}(x)$ is analytic for $|x| < |p|^{-\frac{1}{2}}$. Then the 
choice ${\cal Y}_+ = {\cal Y}_- = f_{N,p,q,M}$ solves the Riemann 
problem (\ref{eq63}) with $|p|^{\frac{1}{2}} < R <
|p|^{-\frac{1}{2}}$.

Exchange relations for the modes of $t(z)$ are then obtained by double
contour integrals.  We shall not give here explicit expressions of the
expansion coefficients of $f_{N,p,q,M}(x)$ but we shall describe 
the connections between the different sector expansions.
\begin{thm}\label{thmeight}
Let $f(x)$ be a meromorphic function of $x^2$ having only simple and 
double poles, denoted by $\pm \alpha_{j}^{-\frac{1}{2}}$, $j\in \NN$.
Suppose that $|\alpha_{1}|^{-\frac{1}{2}}>1$ and 
$|\alpha_{j}|^{-\frac{1}{2}}<|\alpha_{j'}|^{-\frac{1}{2}}$ if 
$j<j'$ and denote $f_l$ the coefficients of the Taylor expansion 
of $f(x)$ for $|x|<|\alpha_{1}|^{-\frac{1}{2}}$:
\begin{equation}
f(x) = \sum_{l=0}^{\infty}f_l\, x^{2l} \, .
\label{eq69}
\end{equation}
Then the relation:
\begin{equation}
\left( \oint_{C_1} \frac{dz}{2\pi iz} \oint_{C_2} 
\frac{dw}{2\pi iw} + \oint_{C_2} \frac{dz}{2\pi iz} \oint_{C_1} 
\frac{dw}{2\pi iw} \right) z^{-n} \, w^{-m} \, [\, f(z/w) \, t(z) \, t(w)
\,-\, f(w/z) \, t(w) \, t(z) \,] = 0 \, ,
\label{eq610}
\end{equation}
where $C_1$ and $C_2$ are circles of radii respectively $R_1$ and
$R_2$ defines a family of exchange relations for the modes $t_n$ 
(\ref{eq326}) of $t(z)$, depending on $R_1/R_2$. \\
If $R_1/R_2 \in [1,|\alpha_{1}|^{\pm \frac{1}{2}}[$ one has:
\begin{equation}
\sum_{l=0}^{\infty}f_l (t_{n-2l}t_{m+2l}-t_{m-2l}t_{n+2l}) = 0 \,,
\label{eq611}
\end{equation}
where $f_l$ is given by (\ref{eq69}), while in the regions 
$R_1/R_2 \in ]|\alpha_{j_0}|^{\pm \frac{1}{2}},|\alpha_{j_0+1}|
^{\pm \frac{1}{2}}[\, $, the exchange relation becomes:
\begin{equation}
\sum_{l\in \ZZ}f^{(j_0)}_l (t_{n-2l}t_{m+2l}-t_{m-2l}t_{n+2l}) = 0 \,,
\label{eq612}
\end{equation}
where $f^{(j_0)}_l$ ($l\in \ZZ$ and $j_0$ here is an index, not a power 
coefficient) is obtained by adding to $f_l$ the contributions:
\begin{equation}
-\sfrac{1}{2} \, \alpha_j^l\, [(1-\alpha_jx^2)f(x)] 
\bigg\vert_{x=\alpha_j^{-1/2}} \,,
\label{eq613}
\end{equation}
for every simple pole $\alpha_j^{-1/2}$ such that $1 \le j \le j_0$, 
and the contributions:
\begin{equation}
-\sfrac{1}{2} \, (l+1)\alpha_j^l\, [(1-\alpha_jx^2)^2f(x)] 
\bigg\vert_{x=\alpha_j^{-1/2}} + \sfrac{1}{4} \, 
\alpha_j^{l-\frac{1}{2}} \,\, \frac{d}{dx}
[(1-\alpha_jx^2)^2 f(x)] \bigg\vert_{x=\alpha_j^{-1/2}} 
\label{eq614}
\end{equation}
for every double pole $\alpha_j^{-1/2}$ such that $1\le j \le j_0$.
\end{thm}
{\bf Proof:} 
If $R_1/R_2 \in [1,|\alpha_1|^{\pm \frac {1}{2}}[\, $, in the integrand 
of (\ref{eq610}) both $|w/z|$ and $|z/w|$ are smaller than $|\alpha
_1|^{-\frac {1}{2}}$: therefore one can develop both $f(z/w)$ and 
$f(w/z)$ in Taylor series according to equation (\ref{eq69}).  Equation 
(\ref{eq611}) then follows immediately.

\medskip

$\bullet$ Suppose now $R_1/R_2 \in \, ]|\alpha_{j_0}|^{\pm \frac{1}{2}},
|\alpha_{j_0+1}|^{\pm \frac{1}{2}}[\, $, where $\alpha_{j_0}^{-1/2}$ 
is a simple pole. Without losing in generality we can suppose  
$R_1/R_2 \in ]|\alpha_{j_0}|^{-\frac{1}{2}},|\alpha_{j_0+1}|^{-\frac{1}{2}}[\,$ 
and rewrite eq. (\ref{eq610}) in the following way:
\begin{equation}
\oint_{C_1} \frac{dw}{2\pi iw} \oint_{C_2} \frac{dz}{2\pi iz} 
z^{-n} \, w^{-m} \, [\, f(z/w) \, t(z) \, t(w) \,-\, f(w/z) \, t(w) \, 
t(z) \,] - (n \leftrightarrow m) = 0 \,.
\label{eq615}
\end{equation} 
In the domain of integration, $f(z/w)$ can be expanded in a convergent 
power series of $z/w$ according to formula (\ref{eq69}), while $f(w/z)$ 
admits the factorization:
\begin{equation}
\displaystyle f \left( \frac{w}{z} \right) =
\prod _{j=1}^{j_0-1} \frac{1}{\displaystyle \left( 1-\alpha_j 
\frac{w^2}{z^2} \right)^{m(j)}} \frac {1}{\displaystyle\left( 1-\alpha_{j_0} 
\frac{w^2}{z^2} \right)} \,\, f^{(j_0)} \left( \frac{w}{z} \right) \,,
\label{eq616}
\end{equation}
where $m(j)=1,2$ is the order of the $j$-th pole and $f^{(j_0)}(w/z)$ is
an analytic function.  The first two factors in the r.h.s. of (\ref{eq616}) 
can be expanded in convergent series of $z/w$:
\begin{equation}
\prod_{j=1}^{j_0-1}\left( 1-\alpha_j \frac{w^2}{z^2} \right)^{-m(j)} 
= \prod_{j=1}^{j_0-1} \,\, \sum_{s=1}^{\infty} g_s^{(j)} \alpha_j^{-s} 
\left( \frac {z}{w} \right)^{2s} \quad , \quad
\left( 1-\alpha_{j_0} \frac{w^2}{z^2} \right)^{-1}
= -\sum_{k=1}^\infty \alpha_{j_0}^{-k} \left( \frac{z}{w} \right)^{2k} 
\label{eq617}
\end{equation}
($g_s^{(j)}=-1$ if $m(j)=1$, $g_s^{(j)}=s-1$ if $m(j)=2$), 
while $f^{(j_0)}(w/z)$ has a convergent series expansion in $w/z$:
\begin{equation}
f^{(j_0)} \left( \frac{w}{z} \right) =
\sum_{l=0}^{\infty} f_l^{(j_0)} \left( \frac{w}{z} \right)^{2l} \,. 
\label{eq618}
\end{equation}
Therefore if $R_1/R_2 \in ]|\alpha_{j_0}|^{\pm
\frac {1}{2}},|\alpha_{j_0+1}|^{\pm \frac {1}{2}}[$ the exchange
relation (\ref{eq610}) can be written, in the sense of formal power series:
\begin{eqnarray}
&& \sum_{l=0}^{\infty}f_l \, t_{n-2l}t_{m+2l} - \oint_{C_1} 
\frac{dw}{2\pi iw} \oint_{C_2} \frac{dz}{2\pi iz} z^{-n} \, w^{-m} \, 
\prod_{j=1}^{j_0-1} \,\, \sum_{s=1}^{\infty} g_s^{(j)} \alpha_j^{-s} 
\left( \frac{z}{w} \right)^{2s} \nonumber \\
&& \times \sum_{k=1}^\infty -\alpha_{j_0}^{-k} \left( \frac{z}{w} \right)^{2k} 
\sum_{l=0}^{\infty} f_l^{(j_0)} \, \left( \frac{w}{z} \right)^{2l} \, t(w) 
\, t(z) - (n \leftrightarrow m) = 0 \,.  
\label{eq619}
\end{eqnarray}
Let us concentrate on the terms:
\begin{equation}
-\sum_{k=1}^\infty \alpha_{j_0}^{-k} \left( \frac{z}{w} \right)^{2k} 
\sum_{l=0}^{\infty}f_l^{(j_0)} \left( \frac{w}{z} \right)^{2l} \,. 
\label{eq620}
\end{equation}
This double sum can be rewritten as a sum of two contributions:
\begin{equation}
-\sum_{h\in \ZZ} \sum_{l=0}^{\infty} f_l^{(j_0)} \alpha_{j_0}^{h-l} 
\left( \frac{w}{z} \right)^{2h} + \sum_{h=0}^{\infty} \sum_{l=0}^h 
f_l^{(j_0)} \alpha_{j_0}^{h-l} \left( \frac{w}{z} \right)^{2h} \,. 
\label{eq621}
\end{equation}
In the second sum of (\ref{eq621}) the term:
\begin{equation}
\sum_{l=0}^h f_l^{(j_0)}\alpha_{j_0}^{h-l} 
\label{eq622}
\end{equation}
is just the coefficient $f_h^{(j_0-1)}$ of the Taylor expansion of the function
$f^{(j_0-1)}(w/z)$ (analytic for $\displaystyle \left\vert\frac {w}{z} 
\right\vert < |\alpha_{j_0}|^{-\frac{1}{2}}$) defined by:
\begin{equation}
f \left( \frac{w}{z} \right) = \prod_{j=1}^{j_0-1}
\frac{1}{\displaystyle\left( 1-\alpha_j \frac{w^2}{z^2} \right)^{m(j)}} 
\,\, f^{(j_0-1)} \left( \frac{w}{z} \right) \,.
\label{eq623}
\end{equation}
Using this result, we can reinsert equation (\ref{eq621})
in (\ref{eq619}). We obtain:
\begin{eqnarray}
&& \sum_{l=0}^{\infty}f_lt_{n-2l}t_{m+2l} - \oint_{C_1} \frac{dw}{2\pi iw} 
\oint_{C_2} \frac{dz}{2\pi iz} z^{-n} \, w^{-m} \, \prod _{j=1}^{j_0-1} 
\,\, \sum_{s=1}^{\infty} g_s^{(j)} \alpha_j^{-s} \left( \frac{z}{w} 
\right)^{2s} \sum_{l=0}^{\infty} f_l^{(j_0-1)} \left( \frac{w}{z} 
\right)^{2l} \, t(w) \, t(z) \nonumber \\
&& + \oint_{C_1} \frac{dw}{2\pi iw} \oint_{C_2} \frac{dz}{2\pi iz} 
z^{-n} \, w^{-m} \, \prod _{j=1}^{j_0-1} \,\, \sum_{s=1}^{\infty} g_s^{(j)} 
\alpha_j^{-s} \left( \frac{z}{w} \right)^{2s} \sum_{h\in \ZZ} 
\sum_{l=0}^{\infty} f_l^{(j_0)} \alpha_{j_0}^{h-l} \left( \frac{w}{z} 
\right)^{2h} \, t(w) \, t(z) \nonumber \\
&& - \, (n \leftrightarrow m) = 0 \,. 
\label{eq624}
\end{eqnarray}
The first line gives the exchange relation for  
$R_1/R_2 \in ]|\alpha_{j_0-1}|^{\pm \frac{1}{2}},|\alpha_{j_0}|
^{\pm \frac{1}{2}}[\, $. The second line is the extra contribution
$(\Delta _{j_0}^1)_{n,m}$ coming from the crossing of the simple singularity 
at $w/z=\alpha_{j_0}^{-1/2}$. The summation over $h$ and $l$ gives:
\begin{eqnarray}
(\Delta_{j_0}^1)_{n,m} &=& \oint_{C_1} \frac{dw}{2\pi iw} \oint_{C_2} 
\frac{dz}{2\pi iz} z^{-n} \, w^{-m} \, \delta \left( \alpha_{j_0} 
\frac{w^2}{z^2} \right) \prod_{j=1}^{j_0-1} \,\, \sum_{s=1}^{\infty} 
g_s^{(j)} \alpha_j^{-s} \left( \frac{z}{w} \right)^{2s} f^{(j_0)} \left( 
\alpha_{j_0}^{-1/2} \right) t(w) \, t(z) \nonumber \\ 
&& - \, (n \leftrightarrow m) \,.
\label{eq625}
\end{eqnarray}
Owing to the properties of the $\delta$-distribution, the series 
expansions in (\ref{eq625}) are all convergent, since 
$|\alpha_{j_0}/\alpha_j|<1$ for $j<j_0$.  Reinserting their sums in 
(\ref{eq625}), one recognizes the expression:
\begin{equation}
\prod _{j=1}^{j_0-1} \frac{1}{\displaystyle \left( 
1-\frac{\alpha_j}{\alpha_{j_0}} \right)^{m(j)}} \,\, f^{(j_0)} \left( 
\alpha_{j_0}^{-1/2} \right) \,,
\label{eq626}
\end{equation}
which, by comparison with (\ref{eq616}), is just the residue of
$f(x)$ at $x=\alpha_{j_0}^{-1/2}$. Then the extra contribution is:
\begin{equation}
(\Delta_{j_0}^1)_{n,m} = 
\oint_{C_1} \frac{dw}{2\pi iw} \oint_{C_2} \frac{dz}{2\pi iz} 
z^{-n} \, w^{-m} \, \delta \left( \alpha_{j_0} \frac{w^2}{z^2} \right) 
[(1-\alpha_{j_0}x^2)f(x)] \Big\vert_{x=\alpha_{j_0}^{-1/2}} 
\,\, t(w) \, t(z) - \, (n\leftrightarrow m) \,,
\label{eq627}  
\end{equation}
or after integration, using the formal power series expansion of
$\delta (x)$:
\begin{equation}
(\Delta_{j_0}^1)_{n,m} = \sum_{l\in \ZZ} \alpha_{j_0}^l 
[(1-\alpha_{j_0}x^2)f(x)] \bigg\vert_{x=\alpha_{j_0}^{-1/2}} 
\,\, (t_{m-2l}t_{n+2l}-t_{n-2l}t_{m+2l}) \,.
\label{eq628}
\end{equation}

\medskip

$\bullet$ Similar but rather longer calculations can be performed in
the case of
a double pole. The extra contribution $(\Delta _{j_0}^2)_{n,m}$ is now:
\begin{eqnarray}
(\Delta _{j_0}^2)_{n,m} &=& 
\oint_{C_1} \frac{dw}{2\pi iw} \oint_{C_2} \frac{dz}{2\pi iz} 
z^{-n} \, w^{-m} \, \delta' \left( \alpha_{j_0} \frac{w^2}{z^2} \right)
[(1-\alpha_{j_0}x^2)^2f(x)] \bigg\vert_{x=\alpha_{j_0}^{-1/2}} 
\,\, t(w) \, t(z) \nonumber \\  
&& -\sfrac{1}{2} \, \oint_{C_1} \frac{dw}{2\pi iw} \oint_{C_2} \frac{dz}{2\pi iz}
z^{-n} \, w^{-m} \, \delta \left( \alpha_{j_0} \frac{w^2}{z^2} \right)
\alpha_{j_0}^{-1/2} \,\, \frac{d}{dx} [(1-\alpha_{j_0}x^2)^2f(x)] 
\bigg\vert_{x=\alpha_{j_0}^{-1/2}} \,\, t(w) \, t(z) \nonumber \\
&& - \, (n\leftrightarrow m) \,, 
\label{eq629}
\end{eqnarray}
where:
\begin{equation}
\delta' (x) = \sum_{h\in \ZZ} h \, x^{h-1} \,.
\label{eq630}
\end{equation}
Integrating (\ref{eq629}) one has:
\begin{eqnarray}
(\Delta _{j_0}^2)_{n,m} &=& \sum_{l\in \ZZ} \left( (l+1)\alpha_{j_0}^l \, 
[(1-\alpha_{j_0}x^2)^2f(x)] \bigg\vert_{x=\alpha_{j_0}^{-1/2}} \right. 
\nonumber \\
&& \left. -\sfrac{1}{2} \, \alpha_{j_0}^{l-\frac{1}{2}} \,\, \frac{d}{dx}
[(1-\alpha_{j_0}x^2)^2f(x)] \bigg\vert_{x=\alpha_{j_0}^{-1/2}} 
\right) (t_{m-2l}t_{n+2l}-t_{n-2l}t_{m+2l}) \,. 
\label{eq631}
\end{eqnarray}
Comparing (\ref{eq631},\ref{eq628}) to the first line of (\ref{eq624}),
one obtains that (\ref{eq613},\ref{eq614}) are precisely the quantities one
has to add to $f_l$ when crossing simple and double singularities.
\finproof

\medskip

{\bf Remark 1:} Theorem \ref{thmeight} gives the mode
exchange algebra between two identical fields $t(z)$ and $t(w)$ using
symmetrized integration contours, see formula (\ref{eq610}).  A similar
result can be obtained for the mode exchange algebra between two
different fields $s_{i}(z)$ and $s_{j}(w)$ but with more complicated
formulae. The Riemann problem (\ref{eq63}) is then generically
solved by two functions ${\cal Y}_{+} \ne {\cal Y}_{-}$; the sector
structure involves poles of {\it both} functions; the separate
contributions of the poles of these two functions to the difference
between adjacent sectors are given by
(\ref{eq613}) and (\ref{eq614}) applied to the corresponding functions.
The mode expansion cannot of course be factored out as in (\ref{eq612})
but takes the generic form
\begin{equation}
\sum_{l=0}^\infty {\cal Y}_{-}(l) \, s_{i}(n-2l) \, 
s_{j}(m+2l) - {\cal Y}_{+}(l) \, s_{j}(m-2l) \, s_{i}(n+2l) = 
0 \,.
\end{equation}
Finally let us emphasize that the symmetrization of the integration contours is
actually not required in the quantum case, even in the case of identical
fields, but it leads to nicer formulae for the coefficients ${\cal
Y}_{\pm}(l)$.

\medskip

{\bf Remark 2:} It may happen that either function ${\cal Y}_{\pm}$,
solution of (\ref{eq63}), have a multiple pole $\alpha$ of order
$\kappa > 2$.  In this case, formulae analogous to (\ref{eq613}) and
(\ref{eq614}) hold: the supplementary contribution to the coefficient
$f_{l}$ is then given by a sum of $\kappa$ terms, each term being
proportional to the $n^{th}$-derivative of $(1-\alpha x^2)^\kappa
f(x)$ taken at $x= \alpha^{-1/2}$, where $n =  0,\dots,\kappa - 1$.

\medskip

{\bf Remark 3:} The condition (\ref{eq68}) is incompatible with the 
classical limit (where $q^{Nh} = p^{1-\beta}$ with $\beta \rightarrow 
0$).  In order to get a classical limit in the case $Nh>2$, it is 
necessary to change (\ref{eq68}) in such a way that the Riemann problem 
(\ref{eq63}) is solved by two functions ${\cal Y}_{+} \ne {\cal 
Y}_{-}$ (in the sense of analytic continuations).

\medskip

{\bf Remark 4:} Finally, let us comment on the classical limit when 
using {\em non-symmetrized} integration contours for the quantum case.  
One starts with the exchange formula between two different fields 
$s_{i}(z)$ and $s_{j}(w)$:
\begin{equation}
{\cal Y}_{-}(z/w) \, s_{i}(z) \, s_{j}(w) = s_{j}(w) \, s_{i}(z) \, 
{\cal Y}_{+}(w/z) \,.
\end{equation}
When using non-symmetrized integration contours $C_{1}$ and $C_{2}$, one 
obtains for the modes
\begin{equation}
\sum_{l=0}^\infty {\cal Y}_{-}(l) \, s_{i}(n-2l) \, 
s_{j}(m+2l) - {\cal Y}_{+}(l) \, s_{j}(m-2l) \, s_{i}(n+2l) = 
0 \,.
\label{eq633}
\end{equation}
Let us set ${\cal Y}_{-}(x) =  1 + \beta f(x)$ and ${\cal Y}_{+}(x) =  1 +
\beta g(x)$. This definition is consistent since the Riemann-Hilbert
problem (\ref{eq63}) reads ${\cal Y}_{+}/{\cal Y}_{-} =  1$ when
$\beta =  0$. (\ref{eq633}) becomes: 
\begin{equation}
\frac{1}{\beta} [s_{i}(n) , s_{j}(m)] + \sum_{l\in\ZZ} f_{l} \, 
s_{i}(n-2l) \, s_{j}(m+2l) - g_{l} \, s_{j}(m-2l) \, 
s_{i}(n+2l) = 0 \,.
\end{equation}
In the classical limit $\beta \rightarrow 0$, one gets therefore, when 
$i \ne j$:
\begin{equation}
\{ s_{i}(n) , s_{j}(m) \} = -\sum_{l\in\ZZ} (f_{l} - g_{-l}) \, 
s_{i}(n-2l) \, s_{j}(m+2l) \,.
\end{equation}
This gives precisely the Poisson structures computed with non-symmetrized 
integration contours with the structure function $f(z/w) - g(w/z)$.  \\
When $i= j$, (\ref{eq633}) can be further decoupled into two exchange relations:
\begin{subequations}
\begin{eqnarray}
&& \sum_{l=0}^\infty ({\cal Y}_{-}(l) + {\cal Y}_{+}(l)) \, 
(s_{i}(n-2l) \, s_{i}(m+2l) -  s_{i}(m-2l) \, s_{i}(n+2l)) = 
0 \,, \\
&& \sum_{l=0}^\infty ({\cal Y}_{-}(l) - {\cal Y}_{+}(l)) \, 
(s_{i}(n-2l) \, s_{i}(m+2l) +  s_{i}(m-2l) \, s_{i}(n+2l)) = 
0 \,.
\end{eqnarray}
\end{subequations}
The first equation leads in the classical case to the Poisson structure
\begin{equation}
\{ s_{i}(n) , s_{i}(m) \} = -\sum_{l\in\ZZ} (f_{l} - g_{-l} 
- f_{-l} + g_{l}) \, s_{i}(n-2l) \, s_{i}(m+2l) \,,
\label{eq638}
\end{equation}
while the second relation has to be interpreted as supplementary 
constraint equations.  (\ref{eq638}) gives precisely the 
antisymmetric Poisson structure obtained in (\ref{eq54}) and 
(\ref{eq56}) by using symmetrized contour integration.  \\
Of course, had we started with symmetrized contours in the quantum case, 
we would get easily the corresponding Poisson structure in the classical 
limit where the mode Poisson bracket is computed with symmetrized 
integration contours.

\section{Conclusions}
\label{sect7}
\setcounter{equation}{0} 

We have studied here exchange relations of the form $f(z/w) \, t(z) \, 
h(w) = h(w) \, t(z) \, f(w/z)$, and their classical Poisson bracket 
limits.  These relations were interpreted in terms of analytic 
continuations and lead us to sets of formal series expansions associated 
to each particular convergent expansion of $f(x)$.

What we have established here is a set of universal structures for
$q$-deformed $W_N$ classical algebras and ${\cal W}_{q,p}(sl(N))$
algebras.These structures still allow for supplementary,
representation-dependent extensions of the type founded in \cite{FF}
or more general!

Nothing was assumed here concerning the fields $t(z)$, $h(w)$, treated
as abstract objects with {\em a priori} no singular behaviour. Of course,
had we constructed explicit realizations of these fields, we would at
one strike in the mode expansion fix a particular sector $(k)$ {\it and} 
add possible central, linear or generically lower-spin extensions as were
constructed in \cite{FF,FR2}.

A closer look at these particular constructions indeed shows that they
all follow from bosonic realizations of the exchange algebra \cite{FR,FF} 
and the extensions are generated either due to cancellations
between otherwise a priori independent fields in the general quadratic
expression \cite{FF}, or (in a different framework, that of deformed
chiral algebras) due to singularities in the operator products
$t(z) \, h(w)$ at the poles of $f(z/w)$, resolved by using explicit
bosonization formulae \cite{FR2}.

It must therefore be expected that extended structures, with lower-spin 
$\delta $-type extensions, play the most important role in practical 
applications of these ${\cal W}_{q,p}(sl(N))$ algebras.  These are also 
crucial to define a reasonable representation theory for these algebras.

Hence the developments from the general frame which we have established 
are twofold.  One should look for explicit realizations (bosonizations ?)  
of the classical/quantum algebras, thereby also getting admissible 
(Jacobi or cocycle solutions) extensions.  Alternatively one should also 
look systematically for allowed extensions of our structures by defining 
and solving the corresponding ``cocycle'' equations.  We hope to report on 
this question in a forthcoming paper \cite {AFRS97c}.

\bigskip

{\bf Acknowledgements}

This work was supported in part by CNRS, Foundation Angelo della Riccia 
and EC network n.  FMRX-CT96-0012.  L.F. is indebted to Centre de 
Recherches Math\'ematiques of Universit\'e de Montr\'eal, where early 
stage of this work was done, for its kind invitation and support.  M.R. 
and J.A. wish to thank ENSLAPP for its kind hospitality.  We also wish 
to thank M. Jimbo and T. Miwa for a number of precious comments and 
indications, and drawing our attention to Ref. \cite{JKOS97}.

\newpage

\renewcommand{\thesection}{\Alph{section}}

\section*{Appendix A: Jacobi theta functions}
\setcounter{section}{1}
\setcounter{equation}{0}

Let $\HH = \{ z\in\CC \,\vert\, \mbox{Im} z > 0 \}$ be the upper
half-plane and $\Lambda_\tau = \{ \lambda_1\,\tau + \lambda_2
\,\vert\, \lambda_1,\lambda_2 \in \ZZ \,, \tau \in \HH \}$ the lattice
with basis $(1,\tau)$ in the complex plane.  One denotes the
congruence ring modulo $N$ by $\ZZ_N \equiv \ZZ/N\ZZ$ with basis
$\{0,1,\dots,N-1\}$.  One sets $\omega = e^{2i\pi/N}$.  Finally, for
any pairs $\gamma=(\gamma_1,\gamma_2)$ and
$\lambda=(\lambda_1,\lambda_2)$ of numbers, we define the
(skew-symmetric) pairing $\langle\gamma,\lambda\rangle \equiv
\gamma_1\lambda_2 - \gamma_2\lambda_1$.

\medskip

\noindent
One defines the Jacobi theta functions with rational characteristics
$\gamma=(\gamma_1,\gamma_2) \in \sfrac{1}{N} \ZZ \times \sfrac{1}{N} \ZZ$ by:
\begin{equation}
\tht\car{\gamma_1}{\gamma_2}(\xi,\tau) = \sum_{m \in \ZZ}
\exp\Big(i\pi(m+\gamma_1)^2\tau + 2i\pi(m+\gamma_1)(\xi+\gamma_2)
\Big) \,.
\label{eqA1}
\end{equation}
The functions $\tht\car{\gamma_1}{\gamma_2}(\xi,\tau)$ satisfy the
following shift properties:
\begin{subequations}
\begin{eqnarray}
&& \tht\car{\gamma_1+\lambda_1}{\gamma_2+\lambda_2}(\xi,\tau) =
\exp(2i\pi \gamma_1\lambda_2) \,\,
\tht\car{\gamma_1}{\gamma_2}(\xi,\tau) \,, \\
&& \tht\car{\gamma_1}{\gamma_2}(\xi+\lambda_1\tau+\lambda_2,\tau) =
\exp(-i\pi\lambda_1^2\tau-2i\pi\lambda_1\xi) \,
\exp(2i\pi\langle\gamma,\lambda\rangle) \,
\tht\car{\gamma_1}{\gamma_2}(\xi,\tau) \,,
\end{eqnarray}
\label{eqA2}
\end{subequations}
where $\gamma=(\gamma_1,\gamma_2) \in \sfrac{1}{N}\ZZ \times
\sfrac{1}{N}\ZZ$ and $\lambda=(\lambda_1,\lambda_2) \in \ZZ \times
\ZZ$.  \\
Moreover, for arbitrary $\lambda=(\lambda_1,\lambda_2)$ (not necessarily
integers), one has the following shift exchange:
\begin{equation}
\tht\car{\gamma_1}{\gamma_2}(\xi+\lambda_1\tau+\lambda_2,\tau) =
\exp(-i\pi\lambda_1^2\tau-2i\pi\lambda_1(\xi+\gamma_2+\lambda_2)) \,
\tht\car{\gamma_1+\lambda_1}{\gamma_2+\lambda_2}(\xi,\tau) \,.
\label{eqA3}
\end{equation}

\medskip

\noindent
Considering the usual Jacobi theta function:
\begin{equation}
\Theta_p(z) = (z;p)_\infty \, (pz^{-1};p)_\infty \, (p;p)_\infty \,,
\label{eqA4}
\end{equation}
where the infinite multiple products are defined by:
\begin{equation}
(z;p_1,\dots,p_m)_\infty = \prod_{n_i \ge 0} (1-z p_1^{n_1}
\dots p_m^{n_m}) \,,
\label{eqA5}
\end{equation}
the Jacobi theta functions with rational characteristics
$(\gamma_1,\gamma_2) \in \sfrac{1}{N} \ZZ \times \sfrac{1}{N} \ZZ$
can be expressed in terms of the $\Theta_{p}$ function as:
\begin{equation}
\tht\car{\gamma_1}{\gamma_2}(\xi,\tau) = (-1)^{2\gamma_1\gamma_2} \,
p^{\frac{1}{2}\gamma_1^2} \, z^{2\gamma_1} \,
\Theta_{p}(-e^{2i\pi\gamma_2} p^{\gamma_1+\frac{1}{2}} z^2) \,,
\label{eqA6}
\end{equation}
where $p = e^{2i\pi\tau}$ and $z = e^{i\pi \xi}$.

\section*{Appendix B: Proof of formula (\ref{eq332})}
\setcounter{section}{2}
\setcounter{equation}{0}

Let us start with the equation (\ref{eq331}):
\begin{equation}
\{ s_{i}(z),s_{j}(w) \} = \sum_{u=-(i-1)/2}^{(i-1)/2} 
\sum_{v=-(j-1)/2}^{(j-1)/2} 
f\Big(q^{v-u}\frac{w}{z}\Big) \, s_{i}(z) \, s_{j}(w) \,.
\end{equation}
where the function $f(x)$ is given by (\ref{eq325}).

\bigskip

\noindent $\bullet$ Let us first suppose that $i<j$.  Then, 
(\ref{eq331}) can be rewritten as:
\begin{equation}
\{ s_{i}(z),s_{j}(w) \} = \sum_{u=1-(i+j)/2}^{(i+j)/2-1} 
\eta(u) f\Big(q^{u}\frac{w}{z}\Big) \, s_{i}(z) \, s_{j}(w) 
\equiv f_{ij}\Big(\frac{w}{z}\Big) \, s_{i}(z) \, s_{j}(w) \,,
\label{eqB2}
\end{equation}
where the function $\eta(u)$ is given by the following graph:
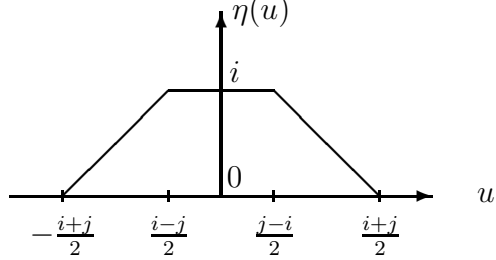
\begin{center}
\begin{picture}(200,70)
\thicklines
\put(-20,0){\vector(1,0){160}}
\put(60,0){\vector(0,1){70}}
\put(0,0){\line(1,1){40}}
\put(40,40){\line(1,0){40}}
\put(80,40){\line(1,-1){40}}
\put(0,-2){\line(0,1){4}}
\put(40,-2){\line(0,1){4}}
\put(80,-2){\line(0,1){4}}
\put(120,-2){\line(0,1){4}}
\put(0,-15){\makebox(0.4,0.6){$-\sfrac{i+j}{2}$}}
\put(40,-15){\makebox(0.4,0.6){$\sfrac{i-j}{2}$}}
\put(80,-15){\makebox(0.4,0.6){$\sfrac{j-i}{2}$}}
\put(120,-15){\makebox(0.4,0.6){$\sfrac{i+j}{2}$}}
\put(160,0){\makebox(0.4,0.6){$u$}}
\put(75,70){\makebox(0.4,0.6){$\eta(u)$}}
\put(65,7){\makebox(0.4,0.6){$0$}}
\put(65,47){\makebox(0.4,0.6){$i$}}
\end{picture}
\end{center}

\vspace{5mm}

\noindent Therefore, one has for the mode Poisson bracket, using 
symmetrized integration contours (with $m,n \in \ZZ$):
\begin{equation}
\{ s_{i}(n),s_{j}(m) \} = \frac{1}{2} \left( \oint_{C_1} 
\frac{dz}{2\pi iz} \oint_{C_2} \frac{dw}{2\pi iw} + \oint_{C_2} 
\frac{dz}{2\pi iz} \oint_{C_1} \frac{dw}{2\pi iw} \right) z^{-n} 
\, w^{-m} f_{ij}\Big(\frac{w}{z}\Big) \, s_{i}(z) \, s_{j}(w) \,.
\label{eqB3}
\end{equation}
One has
\begin{eqnarray}
&&f_{ij}(x) = -2\ln q \sum_{u=1-(i+j)/2}^{(i+j)/2-1} \eta(u) \left\{
\sum_{\ell \ge 0} \left( \frac{2x^2q^{2N\ell+2u}}{1-x^2q^{2N\ell+2u}}
- \frac{x^2q^{2N\ell+2u+2}}{1-x^2q^{2N\ell+2u+2}} \right.
- \frac{x^2q^{2N\ell+2u-2}}{1-x^2q^{2N\ell+2u-2}} \right)
\nonumber \\
&&\quad 
- \sum_{\ell \ge 0} \left( \frac{2x^{-2}q^{2N\ell-2u}}{1-x^{-2}q^{2N\ell-2u}}
- \frac{x^{-2}q^{2N\ell-2u+2}}{1-x^{-2}q^{2N\ell-2u+2}}
- \frac{x^{-2}q^{2N\ell-2u-2}}{1-x^{-2}q^{2N\ell-2u-2}} \right) 
- \frac{x^2q^{2u}}{1-x^2q^{2u}} \nonumber \\
&&\quad \left. + \sfrac{1}{2} \frac{x^2q^{2u+2}}{1-x^2q^{2u+2}}
+ \sfrac{1}{2} \frac{x^2q^{2u-2}}{1-x^2q^{2u-2}} 
+ \frac{x^{-2}q^{-2u}}{1-x^{-2}q^{-2u}} 
- \sfrac{1}{2} \frac{x^{-2}q^{-2u+2}}{1-x^{-2}q^{-2u+2}} 
- \sfrac{1}{2} \frac{x^{-2}q^{-2u-2}}{1-x^{-2}q^{-2u-2}} \right\}
\label{eqB4}
\end{eqnarray}
In the following, one has to distinguish the cases $i+j$ odd and $i+j$ even.
In the sector $k=0$, the integration contours $C_{1}$ and $C_{2}$ of 
radii $R_{1}$ and $R_{2}$ correspond to the choice $\displaystyle 
\frac{R_{1}}{R_{2}} \in \left] \vert q \vert^{1/2},\vert q \vert^{-1/2} 
\right[ \,$ when $i+j$ is odd and $\displaystyle 
\frac{R_{1}}{R_{2}} \in \left] \vert q \vert,\vert q \vert^{-1} 
\right[ \,$ when $i+j$ is even (recall that $\vert q \vert < 1$).  
Hence, in the sector $k=0$, the singularities of the structure function 
lie at $x = q^{\pm u}$, $q^{\pm u \pm 1}$.  

\bigskip

\noindent
\underline{{\bf a) case $i+j$ odd}}  \\
In that case, $u$ takes only half-integer values. One has therefore
\begin{eqnarray}
f_{ij}(x) &=& -2\ln q \left\{ -\sum_{u=1-(i+j)/2}^{(i+j)/2-1} \eta(u) 
\sum_{\ell \ge 1} \sum_{r \ge 1} (q^{r}-q^{-r})^{2} (x^{2r} 
q^{2N\ell r+2ur} - x^{-2r} q^{2N\ell r-2ur}) \right. \nonumber \\ 
&& + \sum_{r \ge 1} x^{2r} \left( \sum_{u>0} 2\eta(u) q^{2ur} - 
\sum_{u>1} \eta(u) q^{2ur-2r} - \sum_{u>-1} \eta(u) q^{2ur+2r} 
\right) \nonumber \\
&& \left.  - \sum_{r \ge 1} x^{-2r} \left( \sum_{u<0} 2\eta(u) 
q^{-2ur} - \sum_{u<1} \eta(u) q^{-2ur+2r} - \sum_{u<-1} \eta(u) 
q^{-2ur-2r} \right) \right\}
\label{eqB5}
\end{eqnarray}
Then resumming the series over $\ell$, (\ref{eqB5}) can be rewritten as:
\begin{eqnarray}
f_{ij}(x) &=& -2\ln q \sum_{r \ge 1} \frac{[r]_{q}}{[Nr]_{q}} 
(q^{r}-q^{-r}) \left\{ x^{2r} \left( q^{Nr} \sum_{u<0} \eta(u) q^{2ur} + 
q^{-Nr} \sum_{u>0} \eta(u) q^{2ur} - \eta(\sfrac{1}{2}) 
\frac{[Nr]_{q}}{[r]_{q}} \right) \right. \nonumber \\
&& \left. - x^{-2r} \left( q^{-Nr} \sum_{u<0} \eta(u) q^{-2ur} + 
q^{Nr} \sum_{u>0} \eta(u) q^{-2ur} - \eta(\sfrac{1}{2}) 
\frac{[Nr]_{q}}{[r]_{q}} \right) \right\}
\label{eqB6}
\end{eqnarray}
{From} the shape of the function $\eta(u)$ (see figure above) it is easy 
to calculate the sums occuring in eq.  (\ref{eqB6}).  One gets
\begin{eqnarray}
&& \sum_{u>0} \eta(u) q^{2ur} = \sum_{u<0} \eta(u) q^{-2ur} =  
\frac{q^{rj}(q^{ri} - q^{-ri}) - i(q^{r}-q^{-r})}{(q^{r}-q^{-r})^{2}}
\nonumber \\
&& \sum_{u>0} \eta(u) q^{-2ur} = \sum_{u<0} \eta(u) q^{2ur} = 
-\frac{q^{-rj}(q^{ri} - q^{-ri}) - i(q^{r}-q^{-r})}{(q^{r}-q^{-r})^{2}}
\label{eqB7}
\end{eqnarray}
Inserting then eq.  (\ref{eqB7}) into eq.  (\ref{eqB6}), one obtains
\begin{equation}
f_{ij}(x) = -2\ln q (q-q^{-1}) \sum_{r \in \ZZ} 
-x^{2r} \frac{[(N-j)r]_{q}[ir]_{q}}{[Nr]_{q}}
\label{eqB8}
\end{equation}
Hence, from eqs. (\ref{eqB3}) and (\ref{eqB8}), formula (\ref{eq332}) 
immediately follows.

\bigskip

\noindent
\underline{{\bf b) case $i+j$ even}} \\
In that case, $u$ takes only integer values; in particular it can take the 
value $u=0$.  One has 
\begin{eqnarray}
f_{ij}(x) &=& -2\ln q \left\{ -\sum_{u=1-(i+j)/2}^{(i+j)/2-1} \eta(u) 
\sum_{\ell \ge 1} \sum_{r \ge 1} (q^{r}-q^{-r})^{2} (x^{2r} 
q^{2N\ell r+2ur} - x^{-2r} q^{2N\ell r-2ur}) \right. \nonumber \\ 
&& + \sum_{r \ge 1} x^{2r} \left( \sum_{u>0} 2\eta(u) q^{2ur} - 
\sum_{u>1} \eta(u) q^{2ur-2r} - \sum_{u>-1} \eta(u) q^{2ur+2r} 
\right) \nonumber \\
&& \left.  - \sum_{r \ge 1} x^{-2r} \left( \sum_{u<0} 2\eta(u) 
q^{-2ur} - \sum_{u<1} \eta(u) q^{-2ur+2r} - \sum_{u<-1} \eta(u) 
q^{-2ur-2r} \right) \right\}
\label{eqB9}
\end{eqnarray}
Although (\ref{eqB9}) has the same form than (\ref{eqB5}), one has to be 
careful in the derivation of this formula due to the peculiar role of 
the value $u = 0$ (in particular note that $i<j$ and $i+j$ even implies 
$\eta(0) = \eta(\pm 1)$).  When resumming the series over $\ell$, one 
obtains after some algebra:
\begin{eqnarray}
f_{ij}(x) &=& -2\ln q \sum_{r \ge 1} \frac{[r]_{q}}{[Nr]_{q}} 
(q^{r}-q^{-r}) \left\{ x^{2r} \left( q^{Nr} \sum_{u<0} \eta(u) q^{2ur} + 
q^{-Nr} \sum_{u>0} \eta(u) q^{2ur} - \eta(0) 
\frac{[(N-1)r]_{q}}{[r]_{q}} \right) \right. \nonumber \\
&& \left. - x^{-2r} \left( q^{-Nr} \sum_{u<0} \eta(u) q^{-2ur} + 
q^{Nr} \sum_{u>0} \eta(u) q^{-2ur} - \eta(0) 
\frac{[(N-1)r]_{q}}{[r]_{q}} \right) \right\}
\label{eqB10}
\end{eqnarray}
The sums over $u$ are now given by:
\begin{eqnarray}
&&\sum_{u>0} \eta(u) q^{2ur} = \sum_{u<0} \eta(u) q^{-2ur} = 
\frac{q^{rj}(q^{ri} - q^{-ri})
- iq^{r}(q^{r}-q^{-r})}{(q^{r}-q^{-r})^{2}} \nonumber \\
&&\sum_{u<0} \eta(u) q^{2ur} = \sum_{u>0} \eta(u) q^{-2ur} = 
-\frac{q^{-rj}(q^{ri} - q^{-ri})
- iq^{-r}(q^{r}-q^{-r})}{(q^{r}-q^{-r})^{2}}
\label{eqB11}
\end{eqnarray}
Inserting (\ref{eqB11}) into (\ref{eqB10}), one gets the same 
formula (\ref{eqB8}) as in the case $i+j$ odd.

\bigskip

\noindent $\bullet$ We consider now the case $i=j$.  Then the 
trapezoidal function $\eta(u)$ degenerates into a triangle:
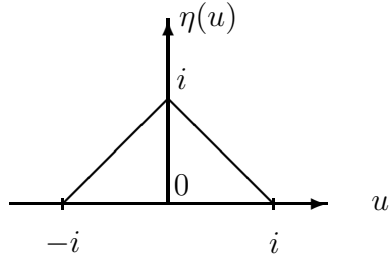
\begin{center}
\begin{picture}(160,70)
\thicklines
\put(-20,0){\vector(1,0){120}}
\put(40,0){\vector(0,1){70}}
\put(0,0){\line(1,1){40}}
\put(40,40){\line(1,-1){40}}
\put(0,-2){\line(0,1){4}}
\put(80,-2){\line(0,1){4}}
\put(0,-15){\makebox(0.4,0.6){$-i$}}
\put(80,-15){\makebox(0.4,0.6){$i$}}
\put(120,0){\makebox(0.4,0.6){$u$}}
\put(55,70){\makebox(0.4,0.6){$\eta(u)$}}
\put(45,7){\makebox(0.4,0.6){$0$}}
\put(45,47){\makebox(0.4,0.6){$i$}}
\end{picture}
\end{center}

\vspace{5mm}

\noindent
The formulae (\ref{eqB2})--(\ref{eqB4}) remain valid with 
the condition $i=j$.
Performing the series expansions in the sector $k=0$, one obtains
\begin{eqnarray}
f_{ii}(x) &=& -2\ln q \left\{ -\sum_{u=1-i}^{i-1} \eta(u) 
\sum_{\ell \ge 1} \sum_{r \ge 1} (q^{r}-q^{-r})^{2} (x^{2r} 
q^{2N\ell r+2ur} - x^{-2r} q^{2N\ell r-2ur}) \right. \nonumber \\ 
&& + \sum_{r \ge 1} x^{2r} \left( \sum_{u>0} 2\eta(u) q^{2ur} - 
\sum_{u>1} \eta(u) q^{2ur-2r} - \sum_{u>-1} \eta(u) q^{2ur+2r} 
+ 1 \right) \nonumber \\
&& \left.  - \sum_{r \ge 1} x^{-2r} \left( \sum_{u<0} 2\eta(u) 
q^{-2ur} - \sum_{u<1} \eta(u) q^{-2ur+2r} - \sum_{u<-1} \eta(u) 
q^{-2ur-2r} + 1 \right) \right\}
\label{eqB12}
\end{eqnarray}
Since now $\eta(0)- \eta(\pm 1) = 1$, eq. (\ref{eqB12}) leads again to 
formula (\ref{eqB10}). Using (\ref{eqB11}), which still holds with 
$i=j$, one gets finally 
\begin{equation}
f_{ii}(x) = -2\ln q (q-q^{-1}) \sum_{r \in \ZZ} 
-x^{2r} \frac{[(N-i)r]_{q}[ir]_{q}}{[Nr]_{q}}
\label{eqB13}
\end{equation}
which achieves the proof of (\ref{eq332}).

\end{document}